\documentclass[leqno,12pt]{amsart}

\usepackage{amsmath,amstext,amssymb,amsopn,amsthm,mathrsfs}
\usepackage{bbm}
\usepackage{verbatim}
\usepackage{enumerate}  
\usepackage{color}
\usepackage{appendix}

\newcommand*\RR{\mathbb{R}}

\newcommand*\NN{\mathbb{N}}
\newcommand*\ZZ{\mathbb{Z}}

\newcommand*\ven{\vert n\vert}
\newcommand*\al{\alpha}
\newcommand*\be{\beta}
\newcommand*\te{\theta}
\newcommand*\va{\varphi}

\newcommand*\fun{\varphi_n}

\newcommand*\jfun{\phi_n^{\al,\be}}
\newcommand*\jfunk{\phi_k^{\al,\be}}

\newcommand*\sjfun{\tilde{\phi}_n^{\al,\be}}
\newcommand*\sjfunk{\tilde{\phi}_k^{\al,\be}}

\newcommand*\jpolk{P_k^{\al,\be}}

\newcommand*\jtpol{\mathcal{P}_n^{\al,\be}}
\newcommand*\jtpolk{\mathcal{P}_k^{\al,\be}}
\newcommand*\jtpoli{\mathcal{P}_{n_i}^{\al_i,\be_i}}

\newcommand*\sjtpolk{\tilde{\mathcal{P}}_k^{\al,\be}}

\newcommand*\tpois{\mathcal{H}_t^{\al,\be}}
\newcommand*\fpois{\mathbb{H}_t^{\al,\be}}

\newcommand*\jRop{\mathbb{R}_r^{\al,\be}}
\newcommand*\jtRop{\mathcal{R}_r^{\al,\be}}

\newcommand*\ve{\varepsilon}
\newcommand*\mf{\mathfrak{m}}

\title[Hardy's inequality for Jacobi expansions]
{Sharp Hardy's inequality for Jacobi and symmetrized Jacobi trigonometric expansions}

\author[P{.} Plewa]{Pawe\l{} Plewa}
\address{Pawe\l{} Plewa \newline
			Faculty of Pure and Applied Mathematics, 
      Wroc\l{}aw University of Science and Technology       \newline
      Wyb{.} Wyspia\'nskiego 27,
      50--370 Wroc\l{}aw, Poland      
      }
\email{pawel.plewa@pwr.edu.pl}

\allowdisplaybreaks

\theoremstyle{plain}
\newtheorem{thm}{Theorem}[section]

\newtheorem{lm}[thm]{Lemma}
\newtheorem{prop}[thm]{Proposition}

\theoremstyle{definition}

\theoremstyle{remark}
\newtheorem*{rem*}{Remark}

\numberwithin{equation}{section}

\theoremstyle{plain}

\usepackage[dvipsnames]{xcolor}

\newcounter{comcount}

\begin{document}
\begin{abstract}
Four Jacobi settings are considered in the context of Hardy's inequality: the trigonometric polynomials and functions, and the corresponding symmetrized systems. In the polynomial cases sharp Hardy's inequality is proved for the type parameters $\al,\be\in(-1,\infty)^d$, whereas in the function systems for $\al,\be\in[-1/2,\infty)^d$. The ranges of these parameters are the widest in which the corresponding orthonormal bases are composed of bounded functions. Moreover, the sharp $L^1$-analogues of Hardy's inequality are obtained with the same restrictions on the parameters $\al$ and $\be$. 
\end{abstract}

\maketitle
\footnotetext{
The paper is a part of author's doctoral thesis written under the supervision of Professor Krzysztof Stempak.
}

\section{Introduction}

In the last two decades many authors investigated Hardy's inequality associated with orthogonal expansions. It states that given an appropriate metric-measure space $(X,\mu)$ and an orthonormal basis $\{\fun\}_{n\in\NN^d}$ in $L^2(X,\mu)$, there holds
\begin{equation}\label{Hardy's_ineq}
\sum_{n\in\NN^d}\frac{|\langle f,\fun\rangle|}{(n_1+\ldots+n_d+1)^E}\lesssim \Vert f\Vert_{H^1(X,\,\mu)},\qquad f\in H^1(X,\mu),
\end{equation}
where $H^1(X,\mu)$ stands for the Hardy space (see Section \ref{Preliminaries}), $\langle\cdot,\cdot\rangle$ denotes the standard inner product in $L^2(X,\mu)$, and $E$ is a positive constant which will be referred to as the admissible exponent.

The inspiration behind \eqref{Hardy's_ineq} was the well known Hardy inequality for the Fourier coefficients (see \cite{HardyLittlewood})
\begin{equation*}
\sum_{k\in\mathbb{Z}} \frac{|\hat{f}(k)|}{|k|+1}\lesssim \Vert f\Vert_{\mathrm{Re} H^1},
\end{equation*}
where $\mathrm{ Re} H^1$ is the real Hardy space. It is formed by the boundary values of the real parts of functions in the Hardy space $H^1(\mathbb{D})$, where $\mathbb{D}$ is the unit disk on the plane.

The study of various inequalities similar to \eqref{Hardy's_ineq} was started by Kanjin \cite{Kanjin1}. He considered two one-dimensional settings: the Hermite functions and the classical Laguerre functions. The obtained admissible exponents were $29/36$ and $1$, respectively. The result for the Hermite expansions was developed by Radha \cite{Radha} and later by Radha and Thangavelu \cite{RadhaThangavelu}. In the latter paper only the multi-dimensional ($d\geq 2$) situation was considered and the received admissible exponent was $E=3d/4$. The case $d=1$ was covered in the article by Z. Li, Y. Yu, and Y. Shi \cite{LiYuShi} with $E=3/4$. In the author's paper \cite{Plewa3} it was proved that those results are sharp, in the sense that the admissible exponent cannot be lower than $3d/4$, $d\geq 1$.

Some analogues of \eqref{Hardy's_ineq} were also considered in the above mentioned settings. Namely, in place of the Hardy space $H^1$ one can consider $H^p$ for $0<p\leq 1$. This was done for the standard Laguerre expansions by Satake \cite{Satake}, for the Hermite functions by Balasubramanian and Radha \cite{BalasRadha}, and by Radha and Thangavelu \cite{RadhaThangavelu} for both systems.

This article is the fourth in a series. In the previous (see \cite{Plewa,Plewa2,Plewa3}) the author investigated the (generalized) Hermite expansions and various Laguerre systems. For information about these settings we refer to \cite{Thangavelu}. In this paper we consider several Jacobi settings. One of them, the trigonometric functions, were considered before by Kanjin and Sato \cite{KanjinSato}, but only for $d=1$.

In the Jacobi systems two type parameters, $\al$ and $\be$, appear. In the multi-dimensional case the admissible range of these parameters in the set $(-1,\infty)^d$. However, for our purposes we require boundedness of the functions constituting the orthonormal bases. Hence, in the function settings, the restriction to $\al,\be\in[-1/2,\infty)^d$ is natural. On the other hand, the Jacobi polynomials belong to $L^\infty$ for all admissible values of the type parameters. In general, dealing with small $\al$ and $\be$, even in the situation $d=1$, namely $\al,\be<-1/2$, is much more complicated because the Jacobi trigonometric polynomials change their asymptotic behaviour (for the details see the proofs of Theorems \ref{Jacobi_pol_nonsym_mainthm} and \ref{nonsym_jpol_L1thm}).

We shall also investigate Hardy's inequality associated with the symmetrized Jacobi trigonometric polynomials and functions. Such expansions were studied in various contexts, see for example the work of Langowski \cite{Langowski}. The procedure of symmetrization was proposed by Nowak and Stempak \cite{NowakStempak(symm)}.

The main results of this paper are sharp Hardy's inequalities, that is Theorems \ref{Jacobi_pol_nonsym_mainthm}, \ref{Jacobi_fun_nonsym_mainthm}, \ref{sym_Jacobi_fun_mainthm}, and  \ref{sym_Jacobi_pol_mainthm}. A method of proving such inequalities in a rather general situations was described in the author's article \cite{Plewa2} (see Theorem \ref{general_thm} below). It is based on some observations made in \cite{LiYuShi}. It consists in estimating kernels of a family of operators $\{R_r\}_{r\in(0,1)}$ (see \eqref{R_r_def} for the definition) directly connected with the considered basis. In the Laguerre and Hermite situations the operators $R_r$ were closely related to the associated heat semigroup. On the other hand, in the Jacobi trigonometric expansions the role is played by the Poisson-Jacobi semigroup.
%It is worth mentioning that in the case of the classical Jacobi polynomials the operators $R_r$ are not so nicely connected with any of the above mentioned semigroups. Fortunately, we are able to deduce the desired results in this setting from the trigonometric one.

The fundamental part of this paper is sharpness of the obtained Hardy's inequalities. It relies on a construction of appropriate counterexamples, which justifies that the admissible exponents cannot be lowered. The essential tool in those proofs are the asymptotic formulas of the Hilb and the Darboux types (see Szeg\"o's monograph \cite{Szego}). Sharpness of Hardy's inequalities was not studied until very recently, when it was done by the author \cite{Plewa2,Plewa3} in the Hermite and Laguerre settings. 

The main obstacle in the analysis of the Jacobi settings is the oscillating nature of these orthonormal bases. In these cases the asymptotic estimates are more complicated than in the Hermite and the Laguerre systems, especially for the type parameters strictly less than $-1/2$. This difficulty is the most conspicuous in the proofs of sharpness of the obtained admissible exponents in Hardy's inequalities. 

Kanjin \cite{Kanjin2} investigated analogues of \eqref{Hardy's_ineq} replacing the $H^1$ norm by the $L^1$ norm. Since $H^1\subset L^1$ the admissible exponent for $L^1$ is always not smaller than the one for $H^1$. In fact, in the cases studied in the author's previous articles, the $L^1$-inequalities hold with any exponents strictly greater than the ones for corresponding Hardy's inequalities (that is excluding the classical Laguerre setting, cf. \cite[Theorem~4.2]{Plewa2}). In this paper we also give the $L^1$-theorems in all of the considered systems.

The article is organized in the following way. In Section \ref{Preliminaries} we introduce the Jacobi systems and give some basic information about the Hardy spaces. Moreover, we describe the method of proving Hardy's inequalities. Section \ref{S3} is devoted to the Jacobi trigonometric polynomials. The sharp inequalities on $H^1$ and $L^1$ are presented. The Jacobi trigonometric function expansions are considered in Section \ref{S4}. Lastly, in Section \ref{S5} we state and prove Hardy's inequalities and theirs $L^1$-analogues in the symmetrized Jacobi trigonometric settings. Additionally, we elaborate on the Hilb and the Darboux type formulas in Appendix.

\subsection*{Notation}
In this article we denote $\NN_+^d=(\NN\setminus\{0\})^d=\{1,2,\ldots\}^d$, where $d\geq 1$ is the dimension. For the real variables we will write $\te,\va$ and $x$ (in both cases $d=1$ and $d\geq 1$). Throughout this paper $k\in\NN$ and $n=(n_1,\ldots,n_d)\in\NN^d$. Moreover, the symbol $\ven$ will stand for the length of the multi-index $n$, namely $\ven=n_1+\ldots+n_d$. The Jacobi type multi-indices $\al,\be>-1$ in the one-dimensional case, whereas $\al=(\al_1,\ldots,\al_d),\be=(\be_1,\ldots,\be_d)\in(-1,\infty)^d$ if $d\geq 1$. The symbol $\lesssim$ is used for inequalities which hold with a multiplicative constant. It may depend on the quantities stated beforehand, but not on the ones quantified afterwards. If both $\lesssim$ and $\gtrsim$ hold simultaneously, then we will denote such relation by $\simeq$. In some estimates we use standard inequalities and write constants using the symbol $C$. It may vary from line to line. If we consider a measure space $(X,\mu)$, where $\mu$ is Lebesgue measure, then we simply write $L^p(X)$ for the Lebesgue spaces and $H^1(X)$ for the Hardy space. 

\subsection*{Acknowledgement}
The author is grateful to Professor Krzysztof Stempak for careful reading the article and his valuable comments. Research supported by the National Science Centre of Poland, NCN grant no. 2018/29/N/ST1/02424.

\section{Preliminaries}\label{Preliminaries}
\subsection{Jacobi settings}
The Jacobi polynomials of order $k\in\NN$ and the type parameters $\al,\be>-1$ are defined on $(-1,1)$ via Rodrigues' formula (see \cite[p.~67]{Szego})
\begin{equation*}
\jpolk(x)\frac{(-1)^k}{2^k k!}(1-x)^{-\al} (1+x)^{-\be} \frac{d^k}{dx^k}\big[ (1-x)^{\al+k} (1+x)^{\be+k}\big].
\end{equation*}
The polynomials $\{\jpolk\}_{k\in\NN}$ form an orthogonal basis in $L^2((-1,1),d\rho_{\al,\be})$, where
\begin{equation*}
d\rho_{\al,\be}(x)=(1-x)^{\al}(1+x)^{\be}\,dx.
\end{equation*}

The normalized Jacobi trigonometric polynomials of order $k\in\NN$ and the type indices $\al,\be>-1$ emerge from the Jacobi polynomials after applying the natural and convenient parametrization $x=\cos \te$, $\te\in(0,\pi)$, and normalization. They are defined by
\begin{equation*}
\jtpolk(\te)= c_k^{\al,\be} \jpolk(\cos \te),\qquad \te\in(0,\pi),
\end{equation*}
with the normalizing constant
\begin{equation}\label{normalizing_constant}
c_k^{\al,\be}=\bigg(\frac{(2k+\al+\be+1)\Gamma(k+\al+\be+1)\Gamma(k+1)}{\Gamma(k+\al+1)\Gamma(k+\be+1)} \bigg)^{1/2},
\end{equation}
where for $k=0$ and $\al+\be=-1$ we write $\Gamma(\al+\be+2)$ in place of $(\al+\be+1)\Gamma(\al+\be+1)$ in the numerator. One can see that $c_k^{\al,\be}\simeq (k+1)^{1/2}$, $k\in\NN$.

We remark that $\jtpolk(\te)=\mathcal{P}_k^{\be,\al}(\pi-\te)$. Analogous symmetries hold also for the other Jacobi systems, which we will introduce below. In the sequel we shall frequently use those identities without any further mention.

Here and later on we define the orthonormal basis only in the one-dimensional situation. However, we can extend those definitions to the case $d\geq 1$ by using the tensor product. For instance, the multi-dimensional Jacobi polynomials of the type parameters $\al,\be\in(-1,\infty)^d$ are given by
\begin{equation*}
\jtpol(\te)=\prod_{i=1}^d \jtpoli (\te_i),\qquad \te=(\te_1,\ldots,\te_d)\in(0,\pi)^d.
\end{equation*}

It is known that for $\al,\be>-1$ (cf. \cite[(7.32.2)]{Szego})
\begin{equation}\label{jtPol_estim}
|\jtpolk(\te)|\lesssim (k+1)^{\max(\al,\be,-1/2)+1/2},\qquad \te\in(0,\pi),\ k\in\NN.
\end{equation}
The family $\{\jtpolk\}_{k\in\NN}$ forms an orthonormal basis in $L^2((0,\pi),d\mu_{\al,\be})$, where
\begin{equation*}
d\mu_{\al,\be}(\te)=\Big(\sin\frac{\te}{2}\Big)^{2\al+1}\Big(\cos\frac{\te}{2}\Big)^{2\be+1}\,d\te.
\end{equation*}

The Jacobi trigonometric polynomials are eigenfunctions of the operator 
\begin{equation*}
\mathcal{J}^{\al,\be}=-\frac{d^2}{d\te^2}-\frac{\al-\be+2\eta\cos\te}{\sin\te}\frac{d}{d\te}+\eta^2,
\end{equation*}
where, here and throughout the entire article, $\eta=(\al+\be+1)/2$. There is
\begin{equation*}
\mathcal{J}^{\al,\be}\jtpolk=\lambda_k^{\al,\be}\jtpolk,
\end{equation*}
where $\lambda_k^{\al,\be}=(k+\eta)^2$.

The operator $\mathcal{J}^{\al,\be}$, considered initially on $C^{\infty}_c(0,\pi)$ (the space of smooth compactly supported functions on $(0,\pi)$), has a natural self-adjoint extension in $L^2((0,\pi),\mu_{\al,\be})$, which we will also denote by $\mathcal{J}^{\al,\be}$. It is given in terms of the eigenfunctions, namely
\[\mathcal{J}^{\al,\be} f = \sum_{k\in\NN}\lambda_k^{\al,\be}\langle f,\jtpolk\rangle \jtpolk,\]
where the associated domain is 
\[\textrm{Dom}\,\mathcal{J}^{\al,\be}=\Big\{f\in L^2((0,\pi),\mu_{\al,\be})\colon \sum_{k\in\NN} \Big(\lambda_k^{\al,\be} |\langle f,\jtpolk\rangle| \Vert\jtpolk\Vert_{L^2((0,\pi),\,\mu_{\al,\be})}\Big)^2<
\infty \Big\}. \]

The associated Poisson-Jacobi semigroup $\{\exp(-t\sqrt{\mathcal{J}^{\al,\be}})\}_{t>0}$ is a semigroup of integral operators. The corresponding kernels are given by
\begin{equation*}
\tpois(\te,\va)=\sum_{k=0}^\infty e^{-t|k+\eta|}\jtpolk(\te)\jtpolk(\va), \qquad \te,\,\va\in(0,\pi),\ t>0.
\end{equation*}
It is known that the kernels $\tpois$ are positive (see for example \cite[p.~189]{NowakSjogrenSzarek}).

A very useful formula for $\tpois$ was obtained by Nowak and Sj\"ogren (see \cite[Proposition~4.1]{NowakSjogren}): for $\al,\be\geq -1/2$ and $\te,\va\in(0,\pi)$ there is
\begin{equation}\label{Pois-Jac_nonsym_pol}
\tpois(\te,\va)=\frac{\sinh \frac{t}{2}}{4^\eta\mu_{\al,\be}(0,\pi)}\iint\limits_{[-1,1]^2} \frac{d\Pi_\al(u)\,d\Pi_\be(v)}{(\cosh \frac{t}{2}-1+q(\te,\va,u,v))^{\al+\be+2}},
\end{equation}
where $\Pi_\al$ and $\Pi_\beta$ are certain measures which we do not need to explicitly define here (for the details see \cite[pp.~189-190]{NowakSjogrenSzarek}), and
\begin{equation*}
q(\te,\va,u,v)=1-u\sin\frac{\te}{2}\sin\frac{\va}{2}-v\cos\frac{\te}{2}\cos\frac{\va}{2},\qquad \te,\va\in(0,\pi),\ u,v\in[-1,1].
\end{equation*}

The integral form of an analogue of $\tpois$ was also obtained in the general case $\al,\be>-1$ in \cite[Theorem~2.1]{NowakSjogrenSzarek}. However, for our purposes much more interesting are the estimates. Such were proved in \cite[Theorem~A.1]{NowakSjogren2} for $\al,\be\geq -1/2$ and then the restraint was removed in \cite[Theorem~6.1]{NowakSjogrenSzarek}. Finally, it was generalized to the derivatives of $\tpois$ (see \cite[Lemma~3.4]{CastroNowakSzarek}) as we state below.

{\thm[Castro, Nowak, Sj\"ogren, Szarek]\label{NowakSjogren_J-P_kernel_estim} If $\al,\be> -1$ and $j\in\NN$, then
\begin{equation*}
\big|\partial_\te^j\tpois(\te,\va)\big|\simeq\big(t^2+\te^2+\va^2\big)^{-\al-1/2} \big(t^2+(\pi-\te)^2+(\pi-\va)^2\big)^{-\be-1/2}\frac{t}{(t^2+(\te-\va)^2)^{1+j/2}},
\end{equation*}
uniformly in $\te,\va\in(0,\pi)$ and $0<t\leq 1$.} 

The next system concerns the Jacobi trigonometric functions. For the type parameters $\al,\be>-1$ the Jacobi functions of order $k\in\NN$ are defined by
\begin{equation}\label{jfunk_definition}
\jfunk(\te)=\Big(\sin\frac{\te}{2}\Big)^{\al+1/2}\Big(\cos\frac{\te}{2}\Big)^{\be+1/2} \jtpolk(\te),\qquad \te\in(0,\pi).
\end{equation}
The system $\{\jfunk\}_{k\in\NN}$ is an orthonormal basis in $L^2((0,\pi))$.

Note that for $\al,\be\geq -1/2$ there is (see \cite[(2.8)]{Muckenhoupt2})
\begin{equation}\label{jfun_estim}
|\jfunk(\te)|\lesssim 1,\qquad k\in\NN,\ \te\in(0,\pi).
\end{equation}

The Jacobi trigonometric functions are the eigenfunctions of the operator
\begin{equation*}
\mathbb{J}^{\al,\be}=-\frac{d^2}{d\te^2}+\frac{\al^2-1/4}{4\sin^2 \frac{\te}{2}}+\frac{\be^2-1/4}{4\cos^2 \frac{\te}{2}},
\end{equation*}
and the corresponding eigenvalues are the same as before, i.e.  $\lambda_k^{\al,\be}=(k+\eta)^2$. This operator, considered initially on $C_c^{\infty}(0,\pi)$, has a natural self-adjoint extension in $L^2(0,\pi)$, given also in terms of the associated eigenfunctions, and is still denoted by $\mathbb{J}^{\al,\be}$. 

Similarly as before, the Poisson-Jacobi semigroup $\{\exp\big(-t\sqrt{\mathbb{J}^{\al,\be}}\big)\}_{t>0}$ is a semigroup of integral operators. The associated kernels are of the form
\begin{equation*}
\fpois(\te,\va)=\sum_{k=0}^{\infty}e^{-t|k+\eta|}\jfunk(\te)\jfunk(\va),\qquad \te,\va\in(0,\pi).
\end{equation*}

Clearly, the following relation between the Poisson-Jacobi kernels in the trigonometric settings holds
\begin{equation}\label{Jacobi_kernel_relation}
\fpois(\te,\va)=\Big(\sin\frac{\te}{2}\sin\frac{\va}{2}\Big)^{\al+1/2}\Big(\cos\frac{\te}{2}\cos\frac{\va}{2}\Big)^{\be+1/2}\tpois(\te,\va),\qquad \te,\va\in(0,\pi).
\end{equation}

\subsection{Symmetrized Jacobi settings}
We shall discuss the symmetrized Jacobi trigonometric systems (cf. \cite{Langowski}). For more information about the symmetrization procedure we refer to \cite{NowakStempak(symm)}.

The one-dimensional symmetrized Jacobi trigonometric polynomials for $k\in\NN$ and $\al,\be>-1$, are defined by
\begin{equation*}
\sjtpolk(\te)=\frac{1}{\sqrt{2}}\left\{ \begin{array}{rl}
\mathcal{P}_{k/2}^{\al,\be}(|\te|),& k\ \text{-- even},\\
\frac{\sin \te}{2}\mathcal{P}_{\lfloor k/2\rfloor}^{\al+1,\be+1}(|\te|), &k\ \text{-- odd},
\end{array}\right.
\end{equation*}
where $\te\in(-\pi,\pi)$ (formally $\te\in(-\pi,0)\cup(0,\pi)$, however for our purposes we may identify it with the full interval). The symbol $\lfloor x\rfloor$ denotes the largest integer less or equal to $x$. The system $\{\sjtpolk\}_{k\in\NN}$ is an orthonormal basis in $L^2((-\pi,\pi),\tilde{\mu}_{\al,\be})$, where
\begin{equation*}
d\tilde{\mu}_{\al,\be}(\te)=\Big| \sin\frac{\te}{2}\Big|^{2\al+1} \Big( \cos\frac{\te}{2}\Big)^{2\be+1}d\te. 
\end{equation*} 

On the other hand, the symmetrized Jacobi trigonometric functions $\{\sjfunk\}_{k\in\NN}$ form an orthonormal basis in $L^2((-\pi,\pi))$ and (for $\in\NN$ and $\al,\be>-1$) are given by
\begin{equation*}
\sjfunk(\te)=\frac{1}{\sqrt{2}}\left\{ \begin{array}{rl}
\phi_{k/2}^{\al,\be}(|\te|),& k\ \text{-- even},\\
\mathrm{sgn}(\te)\phi_{\lfloor k/2\rfloor}^{\al+1,\be+1}(|\te|), &k\ \text{-- odd}.
\end{array}\right.
\end{equation*}

We will not need the explicit formulas for the corresponding Poisson-Jacobi kernels. For the definitions and properties of these systems we refer to \cite{Langowski}.

\subsection{Hardy's inequality}
Throughout this paper we will consider four metric-measure spaces: $(0,\pi)^d$ and $(-\pi,\pi)^d$ both with Lebesgue measure and the measure $\mu_{\al,\be}$ (or $\tilde{\mu}_{\al,\be}$). All of them, equipped with the Euclidean metric, are the spaces of homogeneous type in the sense of Coifman-Weiss (see \cite[pp.~587-588]{CoifmanWeiss}). Hence, we will define the corresponding atomic Hardy spaces accordingly.

Let $(X,\mu)$ be one of the above mentioned metric-measure spaces. A measurable function $a$ supported in a ball $B\subset X$ is called an $H^1(X,\mu)$-atom if $\int a(x)\,d\mu(x)=0$ and $\Vert a\Vert_{L^2(X,\,\mu)}\leq \mu(B)^{-1/2}$. If $\mu(X)<\infty$ we have also an additional atom $a\equiv\mu(X)^{-1}$.

The atomic Hardy space $H^1(X,\mu)$ is composed of all functions $f\in L^1(X,\mu)$ admitting an atomic decomposition, i.e. for a sequence of $H^1(X,\mu)$-atoms $\{a_i\}_{i\in\NN}$ and a sequence of complex numbers $\{\lambda_i\}_{i\in\NN}$, there is
\begin{equation*}
f=\sum_{i\in\NN} \lambda_i a_i,
\end{equation*}
where the series is convergent in $L^1(X,\mu)$. The space $H^1(X,\mu)$ equipped with the norm
\begin{equation*}
\Vert f\Vert_{H^1(X,\,\mu)}=\inf \sum_{i\in\NN} |\lambda_i|,
\end{equation*}
where the infimum is taken over all atomic decompositions of $f$, is a Banach space.

In order to prove Hardy's inequality associated with various Jacobi expansions we shall use \cite[Theorem~2.2]{Plewa2}. For the reader's convenience we state it below.

Let $X$ be an open convex subset of $\RR^d$ and $\mu$ a doubling measure satisfying (lower Ahlfors' condition)
\begin{equation}\label{parameter_N}
\mu(B(x,\rho))\gtrsim \rho^N,
\end{equation}
uniformly in $x\in X$ and $\rho\in (0,\mathrm{diam} X)$, for some $N>0$. Moreover, let $\{\varphi_n\}_{n\in\NN^d}$ be an orthonormal basis in $L^2(X,\mu)$ such that $\varphi_n\in L^\infty (X,\mu)$, $n\in\NN^d$.

We introduce the family of operators $\{R_r\}_{r\in(0,1)}$ defined via 
\begin{equation}\label{R_r_def}
R_r f=\sum_{n\in\NN^d} r^{\ven} \langle f,\fun\rangle \fun,\qquad r\in(0,1),\ f\in L^2(X,\mu),
\end{equation}
where $\langle\cdot,\cdot\rangle$ stands for the standard inner product in $L^2(X,\mu)$. We also impose the assumption that $R_r$ are integral operators and denote the corresponding kernels by $R_r(x,y)$, $x,y\in X$, $r\in(0,1)$.

{\thm\label{general_thm} Let $X,\mu$, $\{\varphi_n\}_{n\in\NN^d}$, and $\{R_r\}_{r\in(0,1)}$ be as above. We assume that there exists $\gamma>0$ and a finite set $\Delta$ constituted of positive numbers satisfying the condition
\begin{equation}\label{parameter_gamma}
\Vert R_r(x,\cdot)-R_r(x',\cdot) \Vert_{L^2(X,\,\mu)}\lesssim \sum_{\delta\in\Delta}|x-x'|^\delta(1-r)^{-\frac{\gamma(N+2\delta)}{N+2}},
\end{equation}
uniformly in $r\in(0,1)$ and almost every $x,x'\in X$ such that $|x'-x|\leq 1/2$.
Then the inequality 
$$\sum\limits_{n\in\NN^d}\frac{|\langle f,\fun\rangle|}{(\ven+1)^E}\lesssim \Vert f\Vert_{H^1(X,\,\mu)} $$
holds uniformly in $f\in H^1(X,\mu)$, where
\begin{equation}\label{E_formula}
E=\frac{\gamma N}{(N+2)}+\frac{d}{2}.
\end{equation}
}

\section{Jacobi polynomial setting}\label{S3}
In this section we establish Hardy's inequality associated with the Jacobi trigonometric polynomials. Unless explicitly stated otherwise, we assume $\al,\be>-1$.

We denote
\begin{equation}\label{jtRop_d=1}
\jtRop(\te,\va)=\sum_{k=0}^{\infty}r^k \jtpolk(\te)\jtpolk(\va),\qquad \te,\va\in(0,\pi),\ r\in(0,1).
\end{equation}
Note that for $\al,\be\geq -1/2$ there  is
\begin{equation}\label{R&H_kernel_relation_pol}
\jtRop(\te,\va)=r^{-\eta} \mathcal{H}^{\al,\be}_{-\log r}(\te,\va).
\end{equation}

Recall that for real non-negative $\omega$ we have
\begin{equation}\label{series_estim_eq}
	\sum_{k=1}^\infty r^k k^\omega\lesssim (1-r)^{-(\omega+1)},\qquad r\in(0,1).
\end{equation}

%\begin{proof}
%	This is known that the for $\omega\in\NN$. Fix a non-integer $\omega>0$. We estimate using H\"older's inequality
%	\begin{align*}
%		\sum_{k=1}^\infty r^k k^{\omega}=\sum_{k=1}^\infty r^{\frac{(\lceil \omega
%		\rceil-\omega)k}{\lceil \omega\rceil}}r^{\frac{k\omega}{\lceil \omega\rceil}}k^\omega\leq \Big(\sum_{k=1}^\infty r^k\Big)^{\frac{\lceil \omega
%		\rceil-\omega}{\lceil \omega \rceil}}\Big(\sum_{k=1}^\infty r^k k^{\lceil \omega
%	\rceil} \Big)^{\frac{\omega}{\lceil \omega\rceil}}\lesssim (1-r)^{-(\omega+1)},
%	\end{align*}
%	uniformly in $r\in(0,1)$. This concludes the proof.
%	 
%\end{proof}

{\prop\label{prop_jpol_R_nonsymmetrized} If $\al,\be> -1$, then
\begin{equation*}
\sup_{\te\in(0,\pi)}\Vert \jtRop(\te,\cdot)\Vert_{L^2((0,\pi),\,\mu_{\al,\be})}\lesssim (1-r)^{-1-\max(\al,\be,-1/2)},
\end{equation*}
and
\begin{equation*}
\sup_{\te\in(0,\pi)}\Big\Vert \frac{\partial}{\partial \te}\jtRop(\te,\cdot)\Big\Vert_{L^2((0,\pi),\,\mu_{\al,\be})}\lesssim(1-r)^{-2-\max(\al,\be,-1/2)},
\end{equation*}
uniformly in $r\in(0,1)$.}
\begin{proof}
Fix $\al,\be>-1$. We invoke Parseval's identity, \eqref{jtPol_estim}, and \eqref{series_estim_eq} obtaining
\begin{align*}
\sup_{\te\in(0,\pi)}\Vert \jtRop(\te,\cdot)\Vert_{L^2((0,\pi),\,\mu_{\al,\be})}&\lesssim \Big(\sum_{k\in\NN} r^{2k} (k+1)^{2\max(\al,\be,-1/2)+1}\Big)^{1/2}\\
&\lesssim (1-r)^{-(\max(\al,\be,-1/2)+1)},\qquad r\in(0,1).
\end{align*}

For the second claim we apply the differentiation formula for $\jtpolk$ (cf. \cite[(4.21.7)]{Szego} or see \cite[(5)]{NowakSjogren})
\begin{equation}\label{jtPol_differentiation_formula}
\frac{d}{d\te}\jtpolk(\te)=-\frac{1}{2}\sqrt{k(k+2\eta)} \mathcal{P}_{k-1}^{\al+1,\be+1}(\te)\sin\te,\qquad k\in\NN,\ \te\in(0,\pi),
\end{equation}
where we set $\mathcal{P}^{\al+1,\be+1}_{-1}\equiv 0$. Hence, in the light of \eqref{jfunk_definition} and \cite[(2.8)]{Muckenhoupt2}
\begin{equation}\label{jtpol_differ_estim}
\Big\vert\frac{d}{d\te}\jtpolk(\te)\Big\vert\lesssim (k+1)^{\max(\al,\be,-1/2)+3/2},\qquad \te\in(0,\pi),\ k\in\NN. 
\end{equation}

Therefore,
\begin{equation*}
\sup_{\te\in(0,\pi)}\Big\Vert \frac{\partial}{\partial \te}\jtRop(\te,\cdot)\Big\Vert_{L^2((0,\pi),\,\mu_{\al,\be})}\lesssim (1-r)^{-(\max(\al,\be)+2)},\qquad r\in(0,1),
\end{equation*}
and this finishes the proof of the proposition.
\end{proof}

{\thm\label{Jacobi_pol_nonsym_mainthm} Let $\al,\be\in(-1,\infty)^d$. Hardy's inequality for the Jacobi trigonometric polynomials holds with the admissible exponent 
$$E=E_{d,\al,\be}=\frac{3d}{2}+\sum_{i=1}^d \max(\al_i,\be_i,-1/2),$$
namely 
\begin{equation*}
\sum_{n\in\NN^d} \frac{|\langle f,\jtpol\rangle|}{(\ven+1)^{E}}\lesssim \Vert f\Vert_{H^1((0,\pi)^d,\,\mu_{\al,\be})},
\end{equation*}
uniformly in $f\in H^1((0,\pi)^d,\,\mu_{\al,\be})$. The admissible exponent is sharp, i.e. for any $\ve>0$ there exists $f\in H^1((0,\pi)^d,\,\mu_{\al,\be})$ such that
\begin{equation*}
\sum_{n\in\NN^d} \frac{|\langle f,\jtpol\rangle|}{(\ven+1)^{E-\ve}}=\infty.
\end{equation*} 
}

\begin{proof}
Fix $\al,\be\in( -1,\infty)^d$. In order to prove the first claim we consider the multi-dimensional version of \eqref{jtRop_d=1}, namely
\begin{equation*}
\jtRop(\te,\va)=\sum_{n\in\NN^d}r^{\ven} \jtpol(\te)\jtpol(\va),\qquad \te,\va\in(0,\pi)^d,\ r\in(0,1).
\end{equation*}
We will invoke Theorem \ref{general_thm}. The proper assumptions are satisfied. Indeed, it can be easily verified, that for the measure $\mu_{\al,\be}$ there is (see \eqref{parameter_N})
\begin{equation*}
	N=2d+2\sum_{i=1}^d \max(\al_i,\be_i,-1/2). 
\end{equation*}
Moreover, Proposition \ref{prop_jpol_R_nonsymmetrized} yields
\begin{align*}
\sup_{\te\in(0,\pi)^d}\Big\Vert\big\vert\nabla_\te \jtRop(\te,\cdot)\big\vert\Big\Vert_{L^2((0,\pi)^d,\,\mu_{\al,\be})}\leq&  \sum_{i=1}^d \sup_{\te_i\in(0,\pi)}\Big\Vert\partial_{\te_i} \mathcal{R}^{\al_i,\be_i}(\te_i,\cdot)\Big\Vert_{L^2((0,\pi),\,\mu_{\al_i,\be_i})}\\
&\times \prod_{\stackrel{j=1}{j\neq i}}^d \sup_{\te_j\in(0,\pi)}\Big\Vert \mathcal{R}^{\al_j,\be_j}(\te_j,\cdot)\Big\Vert_{L^2((0,\pi),\,\mu_{\al_j,\be_j})}\\
\lesssim &(1-r)^{-(d+1)-\sum_{i=1}^d\max(\al_i,\be_i,-1/2)},\qquad r\in(0,1).
\end{align*}
Hence, the appropriate version of \eqref{parameter_gamma} holds with $N$ as above, $\Delta=\{1\}$, and
$$\gamma=d+1+\sum_{i=1}^d\max(\al_i,\be_i,-1/2).$$
 
Therefore, by \eqref{E_formula} we receive
$$E=\frac{3d}{2}+\sum_{i=1}^d \max(\al_i,\be_i,-1/2).$$

Now we pass to the essential part of the proof, namely we shall justify sharpness. In the one-dimensional case, we will construct, for a given $K\in\NN_+$, an appropriate $H^1((0,\pi),\mu_{\al,\be})$-atom $a_K$ such that
\begin{equation}\label{H1_jtpol_claim}
\sum_{k=1}^K \frac{|\langle a_K,\jtpolk\rangle|}{k^{3/2+\max(\al,\be,-1/2)-\ve}}\gtrsim K^{\ve/2},
\end{equation}
where the implicit multiplicative constant does not depend on $K$. This is indeed sufficient; one may apply for instance the uniform boundedness principle. 

Without any loss of generality we assume that $\al\geq \be>-1$. Firstly, we consider the case $\al\geq -1/2$. For fixed large $K\in\NN$, sufficiently small $c>0$, and $0<\delta<1$, we set
\begin{equation*}
a_K(\te)=\left\{ \begin{array}{rl}
-(K/c)^{2\al+2}(1-\delta^{2\al+2}),&\te\in(0,c\delta K^{-1}],\\
C_{\delta,K}(K\delta/c)^{2\al+2}, &\te\in(c\delta K^{-1},cK^{-1}],\\
0,&\text{otherwise},
\end{array}
\right.
\end{equation*}
where $C_{\delta,K}$ is uniquely chosen in such way that $\int a(\te)\,d\mu_{\al,\be}(\te)=0$. Observe that though $C_{\delta,K}$ depends on $K$, we have the estimate $C_{\delta,K}\simeq 1$, $K\in\NN_+$, where the implicit constants depend only on $\al$, $\be$, and $\delta$. It is easy to check that for sufficiently small $\delta$ (depending on $\al$ and $\be$) the function $a_K$ is an $H^1((0,\pi),\mu_{\al,\be})$-atom.

Note that \eqref{eq:A3} implies that $\frac{d}{d\te}\jtpolk(\te)$ is negative for sufficiently small $\te$. Hence, invoking the mean value theorem and \eqref{eq:A4}, we receive for $K/2\leq k\leq K$ and sufficiently large $K$ the estimate
\begin{align*}
\Big| \int_0^{\frac{c}{K}}a_K(\te)\jtpolk(\te)\, d\mu_{\al,\be}(\te)\Big|&=\Big| \int_0^{\frac{c}{K}}a_K(\te)\big(\jtpolk(\te)-\jtpolk(c\delta K^{-1})\big)\, d\mu_{\al,\be}(\te)\Big|\\
&=\int_0^{\frac{c}{K}}|a_K(\te)|\Big(-\frac{d}{d\te}\jtpolk(\xi_\te)\Big) |\te-c\delta K^{-1}|\, d\mu_{\al,\be}(\te)\\
&\geq B (1-\delta^{2\al+2}) k^{\al+5/2}\Big(\frac{K}{c}\Big)^{2\al+2}\int_0^{\frac{c\delta }{K}}\big(\frac{c\delta}{K}-\te\big)\te^{2\al+2}\, d\te\\
&\qquad-Ck^{\al+1/2} K^{-1} \Vert a\Vert_{L^1((0,\pi),\,\mu_{\al,\be})}\\
&\gtrsim k^{\al+5/2} K^{-2}\\
&\gtrsim K^{\al+1/2},
\end{align*}
where, for $\te\in(0,c/K)$, $\xi_\te$ lies between $\te$ and $c\delta K^{-1}$, $B$ is the constant from \eqref{eq:A4}, and $C$ is the constant emerging from the estimate on the reminder appearing in \eqref{eq:A3}. Thus,
\begin{equation*}
\sum_{k=1}^{K} \frac{|\langle a_K,\jtpolk\rangle|}{k^{\frac{3}{2}+\al-\ve}}\gtrsim K^{\al+1/2} \sum_{k=\lceil K/2\rceil}^{K} k^{-3/2-\al+\ve}\simeq K^\ve,\qquad K\in\NN_+,
\end{equation*}
where $\lceil x\rceil$ denotes the smallest integer not less than $x$. This concludes the proof of \eqref{H1_jtpol_claim} for $\max(\al,\be)\geq -1/2$.

Now, let us assume that $-1/2> \al\geq \be>-1$. We define for a given $K\in\NN_+$ and $c$ as above
\begin{equation*}
b_K(\te)=\left\{
\begin{array}{rl}
-C_KK^{2\al+2-\ve/2},&\te\in(0,cK^{-1}],\\
K^{1-\ve/2},&\te\in(\pi/2-c K^{-1},\pi/2+cK^{-1}],\\
0,& \text{otherwise}.
\end{array}\right.
\end{equation*}
This time, similarly as before, $C_K$ is uniquely chosen in such way that $\int_0^\pi  b_K(\te) \,d\mu_{\al,\be}(\te)=0$. Moreover, $C_K\simeq 1$, $K\in\NN_+$. Note that $\mathrm{supp}\,b_K\subset [0,\pi/2+cK^{-1}]$ and $\mu_{\al,\be}([0,\pi/2+cK^{-1}])\simeq 1,$ $K\in\NN_+$. One may check that $\Vert b_K\Vert_{L^q((0,\pi),\,\mu_{\al,\be})}\simeq K^{1-1/q-\ve/2}$, $K\in\NN_+$. Hence, for $1<q<\frac{2}{2-\ve}$ and sufficiently large $K$, we have 
\[\Vert b_K\Vert_{L^q((0,\pi),\,\mu_{\al,\be})}\leq  \big(\mu_{\al,\be}([0,\pi/2+cK^{-1}])\big)^{1/q-1}.\]
Fix $q$ satisfying the above-mentioned conditions. Observe that $b_K$ is an $(1,q)$-atom (cf. \cite[p.~591]{CoifmanWeiss}). In the light of \cite[Theorem~A]{CoifmanWeiss} in order to verify sharpness it suffices to justify the appropriate version of \eqref{H1_jtpol_claim} for $(1,q)$-atom, namely
\begin{equation}\label{thm:jtpolk:claim}
\sum_{k=1}^{K}\frac{|\langle b_K,\jtpolk\rangle|}{k^{1-\ve}}\gtrsim K^{\ve/2},\qquad K\in\NN_+.
\end{equation}

%We shall use the Darboux formula (cf. \cite[(8.21.18)]{Szego})
%\begin{equation}\label{Darboux_formula}
%\jtpolk(\te)\simeq  \cos(k\te+\eta)+ r^{\al,\be}_k,\qquad \te\in[\pi/8,7\pi/8],
%\end{equation}
%where $\eta=\pi(\be-\al)/4$ and $|r^{\al,\be}_k|\lesssim (k+1)^{-1}$.

We will apply \eqref{eq:A5}. Note that under our assumptions $\eta\in(-1/2,0)$ and $(2\al+1)\pi/4\in(-\pi/4,0)$. Moreover, we see that for $\te\in (\pi/2-c K^{-1},\pi/2+cK^{-1})$ and $k=4l\leq K$, $l\in\NN_+$, we have
\begin{equation*}
(4l+\eta)\te-\frac{(2\al+1)\pi}{4}\in\Big(2\pi l-\frac{\pi}{4}-c,2\pi l+\frac{\pi}{4}+c \Big)
\end{equation*}
and thus
\begin{equation*}
	\cos\Big((4l+\eta)\te-\frac{(2\al+1)\pi}{4}\Big)\gtrsim 1,\qquad l\leq \lfloor K/4\rfloor,\ \te\in\Big(\frac{\pi}{2}-\frac{c}{K},\frac{\pi}{2}+\frac{c}{K}\Big)
\end{equation*}
provided that $c$ is sufficiently small.

Hence, invoking the above and \eqref{eq:A2,5} along with \eqref{eq:A2} we have for $K/2\leq k=4l\leq K$, $l\in\NN_+$, the estimate
\begin{align*}
\Big| \int_0^{\pi}b_K(\te)\jtpolk(\te)\, d\mu_{\al,\be}(\te)\Big|&\geq\sqrt{\frac{2}{\pi}} K^{1-\ve/2}\int_{\pi/2-c K^{-1}}^{\pi/2+cK^{-1}}\cos\Big((4l+\eta)\te-\frac{(2\al+1)\pi}{4}\Big)\\
&\hspace{1cm} \times \Big(\sin\frac{\te}{2}\Big)^{-\al-1/2}\Big(\cos\frac{\te}{2}\Big)^{-\be-1/2}\, d\mu_{\al,\be}(\te)\\
&\hspace{1cm} -\frac{C_K}{A} K^{2\al+2-\ve/2}k^{\al+1/2}\int_0^{\frac{c}{K}}\, d\mu_{\al,\be}(\te)\\
&\hspace{1cm} -C(k^{-1}+k^{\al-3/2})\Vert b\Vert_{L^1((0,\pi),\,\mu_{\al,\be})}\\
&\gtrsim K^{-\ve/2} -K^{-\ve/2}k^{\al+1/2}\\
&\gtrsim K^{-\ve/2},
\end{align*}
uniformly in sufficiently large $K\in\NN_+$, where $A$ is the constant from \eqref{eq:A2}, and $C$ emerges from estimates on the remainders in \eqref{eq:A2,5} and \eqref{eq:A5}. Note that in the last inequality we used the fact that $\al+1/2<0$. Thus,
\begin{equation*}
\sum_{k=1}^K \frac{|\langle b_K,\jtpolk\rangle|}{k^{1-\ve}}\gtrsim K^{-\ve/2}\sum_{l=\lceil K/8\rceil}^{\lfloor K/4\rfloor} (4l)^{-1+\ve} \gtrsim K^{\ve/2},\qquad K\in\NN_+.
\end{equation*}
This concludes verification of \eqref{thm:jtpolk:claim} and, consequently, finishes the proof of sharpness for $d=1$.

In the multi-dimensional case we may also assume that $\al_i\geq \be_i$ for $i=1,\ldots,d$. The appropriate counterexample of sequence of atoms is defined by the tensor product of the above atoms, namely the function
\[\mathbf{a}_K(x)=\prod_{i=1}^d (\mathbf{a}_K)_i(x_i),\]
where $(\mathbf{a}_K)_i(x_i)=a_K(x_i)$ if $\al_i\geq -1/2$ and $(\mathbf{a}_K)_i(x_i)=b_K(x_i)$ if $\al_i< -1/2$, with the slight changes:
\begin{itemize}
	\item in place of the powers $2\al+2-\ve/2$ and $1-\ve/2$ of $K$ in the definition of $b_K$ we put $2\al+2-\ve/(d+1)$ and $1-\ve/(d+1)$, respectively;
	\item  in place of the powers $2\al+2$ of $K/c$ and of $K\delta/c$ in the definition of $a_K$ we put $2\al+2-\ve/(d+1)$.
\end{itemize}
Then we observe that
\[ \Vert \mathbf{a}_K\Vert_{L^q((0,\pi)^d,\,\mu_{\al,\be})}\simeq K^{(1-1/q)\sum_{i=1}^d \max(1,2+2\al_i)-\frac{d\ve}{d+1}},  \]
and $\mu_{\al,\be}(\mathrm{supp}\,\mathbf{a}_K)\lesssim 1$. Hence, $\mathbf{a}_K$ is an $(1,q)$-atom for 
\[q\in(1,2]\cap\bigg(1,\frac{(d+1)\sum_{i=1}^d\max(1,2+2\al_i) }{(d+1)\sum_{i=1}^d \max(1,2+2\al_i)-d\ve}\bigg). \]
We omit the details.

%Let $j=|\{i\colon \al_i<-1/2\}|$. If $j=0$, then it suffices to justify
%\begin{equation*}
%\sum_{n\in\NN^d} \frac{|\langle \mathbf{a},\jtpol\rangle|}{(\ven+1)^{3d/2+|\al|-\ve}}\gtrsim K^{\ve/(d+1)},
%\end{equation*}
%where $|\al|=\al_1+\ldots+\al_d$. Indeed, we have
%\begin{align*}
%	\sum_{n\in\NN^d} \frac{|\langle \mathbf{a},\jtpol\rangle|}{(\ven+1)^{3d/2+|\al|-\ve}}\gtrsim K^{-2d}\sum_{n\in\{1,\ldots,K\}^d} \frac{\prod_{i=1}^d (n_i+1)^{\al_i+5/2}}{(\ven+1)^{3d/2+|\al|-\ve}}
%\end{align*}
%
%Note that $\mathbf{a}$ is an $(1,p)$-atom, where
%\[p=\min\Big(\frac{2\al_i+2}{2\al_i+2-\ve/(d+1)}\colon i=1,\ldots,d\Big)\].
%
%Hence,
%\begin{align*}
%	\sum_{n\in\NN^d} \frac{|\langle f,\jtpol\rangle|}{(\ven+1)^{3d/2+|\al|-\ve}} 
%\end{align*}
\end{proof}

Before we will present an $L^1$-analogue of Theorem \ref{Jacobi_pol_nonsym_mainthm}, we give an auxiliary lemma.

{\lm\label{series_divergence_lemma} Let $d\geq 1$ and $j\in\{0,1,\ldots,d\}$. Assume that $a_i,b_i\in\RR$, $i=1,\ldots,j,$ and $a_i$ are such that $a_i\neq l\pi$, $l\in\ZZ$. Then for any $\omega_s\in\RR$, $s=j+1,\ldots,d$, there is
\begin{equation}\label{series_divergent_claim1}
\sum_{n\in\NN^d} \frac{\prod_{i=1}^j |\cos(a_in_i+b_i)| \prod_{s=j+1}^d (n_s+1)^{\omega_s}}{(\ven+1)^{d+\sum_{s=j+1}^d\omega_s}}=\infty,
\end{equation}
If we replace (all or some of) the cosines by the sines, then the claim holds as well.}
\begin{proof}
Clearly, if $j=0$ ($j=d$), then we replace the first (second) product in the numerator by $1$. We shall prove the claim using the induction over $d$.

%For the sake of clarity we can assume that $a_i=1$, $b_i=0$, for $i=1,\ldots,j$.

If $d=1$, then we have to consider two cases:
\begin{equation*}
\sum_{k\in\NN} \frac{|\cos (ak+b)|}{k+1},\qquad \sum_{k\in\NN}\frac{1}{k+1},
\end{equation*}
for some $a\neq l\pi$, $l\in\ZZ$, and $b\in\RR$. Clearly, both series above are divergent. 

In order to perform the inductive  step we fix $d\geq 1$ and assume that for $1,\ldots,d$, each of the series from \eqref{series_divergent_claim1} diverges, and we will show the corresponding property for $d+1$ (all under the appropriate assumption on $a_i$'s). Fix $j\in\{0,\ldots,d+1\}$ and $\omega_s\in\RR$ for $s=j+1,\ldots,d+1$. Let $w=|\{s\colon \omega_s<0\}|$.

If $w\geq 1$, then we may assume that $\omega_{d+2-w},\ldots,\omega_{d+1}<0$. For the time being we denote $\ven=n_1+\ldots+n_{d+1}$. Observe that
\begin{align*}
\sum_{n_{d+2-w},\ldots,n_{d+1}\in\NN} \frac{\prod_{l=d+2-w}^{d+1}(n_l+1)^{\omega_l}}{(1+\ven)^{d+1+\sum_{s=j+1}^{d+1}\omega_s}}&\geq \sum_{n_{d+2-w},\ldots,n_{d+1}\in\NN}(1+\ven)^{-d-1-\sum_{s=j+1}^{d+1-w}\omega_s}\\
&\gtrsim \big(n_1+\ldots+n_{d+1-w}+1\big)^{-d-1+w-\sum_{s=j+1}^{d+1-w}\omega_s},
\end{align*}
where in the second inequality we applied the simple fact that
\begin{equation}\label{series_simple_estim}
	\sum_{k\in\NN} (k+K)^{-\tau}\geq \frac{1}{\tau-1} K^{-\tau+1},
\end{equation}
for $K>0$ and $\tau>1$. Here and later on we use the convention that $\sum_{s=i_1}^{i_2}t_s=0$ for $i_1>i_2$ and any sequence $\{t_s\}_{s\in\NN}$. Hence,
\begin{align*}
&\sum_{n\in\NN^{d+1}} \frac{\prod_{i=1}^j |\cos (a_in_i+b_i)| \prod_{s=j+1}^{d+1} (n_s+1)^{\omega_s}}{(\ven+1)^{d+1+\sum_{s=j+1}^{d+1}\omega_s}}\\
&\hspace{1cm}\gtrsim \sum_{n_1,\ldots,n_{d+1-w}\in\NN} \frac{\prod_{i=1}^j |\cos (a_in_i+b_i)| \prod_{s=j+1}^{d+1-w} (n_s+1)^{\omega_s}}{(n_1+\ldots+n_{d+1-w}+1)^{d+1-w+\sum_{s=j+1}^{d+1-w}\omega_s}}.
\end{align*}
The inductive hypothesise yields the claim. 

Now consider the situation $w=0$. For the time being denote $n=(n',n_{d+1})\in\NN^d\times \NN$, and $|n'|=\ven-n_{d+1}$. We distinguish two cases: $j=0$ or $j\geq 1$. In the former we calculate using \eqref{series_simple_estim}
\begin{align*}
\sum_{n\in\NN^{d+1}}& \frac{ \prod_{s=1}^{d+1} (n_s+1)^{\omega_s}}{(\ven+1)^{d+1+\sum_{s=1}^{d+1}\omega_s}}\\
&\geq \sum_{n'\in\NN^d} \Big(\prod_{s=1}^d (n_s+1)^{\omega_s}\Big) \sum_{n_{d+1}=\lfloor |n'|/d\rfloor}^\infty \Big(\frac{n_{d+1}+1}{\ven+1}\Big)^{\omega_{d+1}}(\ven+1)^{-d-1-\sum_{s=1}^d \omega_s}\\
&\gtrsim \sum_{n'\in\NN^d} \Big(\prod_{s=1}^d (n_s+1)^{\omega_s}\Big) \sum_{n_{d+1}\in\NN} \Big(n_{d+1}+\frac{1}{2}|n'|+2\Big)^{-d-1-\sum_{s=1}^d \omega_s}\\
&\gtrsim \sum_{n'\in\NN^d} \frac{ \prod_{s=1}^{d} (n_s+1)^{\omega_s}}{(|n'|+1)^{d+\sum_{s=1}^{d}\omega_s}},
\end{align*}
and the latter series is divergent due to the inductive hypothesis. For the second case, $j\geq 1$, we first remark that 
\begin{align*}
|\cos (a_1n_1+b_1)|\geq \cos^2 (a_1 n_1+b_1)=\frac{1+\cos (2a_1 n_1+2b_1)}{2}.
\end{align*}
Thus,
\begin{align*}
\sum_{n\in\NN^{d+1}}&  \frac{\prod_{i=1}^j |\cos (a_i n_i+b_i)| \prod_{s=j+1}^{d+1} (n_s+1)^{\omega_s}}{(\ven+1)^{d+1+\sum_{s=j+1}^{d+1}\omega_s}}\\
&\geq \sum_{n\in\NN^{d+1}}\frac{\prod_{i=2}^j |\cos (a_i n_i+b_i)| \prod_{s=j+1}^{d+1} (n_s+1)^{\omega_s}}{2(\ven+1)^{d+1+\sum_{s=j+1}^{d+1}\omega_s}}\\
&\qquad+\sum_{n\in\NN^{d+1}} \frac{\cos (2a_1 n_1+2b_1) \prod_{i=2}^j |\cos (a_i n_i+b_i)| \prod_{s=j+1}^{d+1} (n_s+1)^{\omega_s}}{2(\ven+1)^{d+1+\sum_{s=j+1}^{d+1}\omega_s}}\\
&\geq \frac{1}{d+\sum_{s=j+1}^{d+1}\omega_s} \sum_{n_2,\ldots,n_{d+1}\in\NN}\frac{\prod_{i=2}^j |\cos (a_i n_i+b_i)| \prod_{s=j+1}^{d+1} (n_s+1)^{\omega_s}}{2(n_2+\ldots+n_{d+1}+1)^{d+\sum_{s=j+1}^{d+1}\omega_s}}\\
&\qquad+\sum_{n_1\in\NN} \cos (2a_1 n_1+2b_1) \sum_{n_2,\ldots,n_{d+1}\in\NN} \frac{ \prod_{i=2}^j |\cos (a_i n_i+b_i)| \prod_{s=j+1}^{d+1} (n_s+1)^{\omega_s}}{2(\ven+1)^{d+1+\sum_{s=j+1}^{d+1}\omega_s}},
\end{align*}
where in the last inequality we again used \eqref{series_simple_estim}. The inductive hypothesis yields that the first series on the right hand side of the latter inequality is divergent. On the other hand, the second one is conditionally convergent (unless $a_1$ and $b_1$ are such that $\cos(2a_1n_1+2b_1)\equiv 0$, $n_1\in\NN$; then it is simply equal to zero). To justify this one may use Dirichlet's test as $\omega_s\geq0$ for all $s$.

We emphasise that if we put sines in place of the cosines, then the estimates are also correct. This concludes the proof of the lemma.
\end{proof}

{\thm\label{nonsym_jpol_L1thm} Let $\al,\be\in(-1,\infty)^d$ and fix $\ve>0$. The following inequality holds
\begin{equation*}
\sum_{n\in\NN^d} \frac{|\langle f,\jtpol\rangle|}{(\ven+1)^{3d/2 +\sum_{i=1}^d \max(\al_i,\be_i,-1/2)+\ve}}\lesssim \Vert f\Vert_{L^1((0,\pi)^d,\,\mu_{\al,\be})},
\end{equation*}
uniformly in $L^1((0,\pi)^d,\mu_{\al,\be}).$ Moreover, the admissible exponent is sharp, i.e. the estimate is not valid for $\ve=0$.}

\begin{proof}
The inequality follows from \eqref{jtPol_estim}. In order to prove sharpness it is sufficient to justify that (see \cite[Lemma~1]{Kanjin2})
\begin{equation}\label{L1_jpol_claim}
\sup_{\te\in(0,\pi)^d}\sum_{n\in\NN_+^d}\frac{|\jtpol(\te)|}{\ven^{3d/2 +\sum_{i=1}^d \max(\al_i,\be_i,-1/2)}}=\infty.
\end{equation}

We begin with the one-dimensional situation. Without loss of generality we assume that $\al\geq\be> -1$. Firstly, we consider the case $\al\geq -1/2$. We apply \eqref{eq:A2,5} and \eqref{eq:A2} for $\te=\frac{c}{K}$ and $k\leq K$, where $K$ is large enough, and get
\begin{equation*}
\sum_{k\in\NN_+} \frac{\big|\jtpolk(\frac{c}{K})\big|}{k^{3/2 +\al}}\geq  A'\sum_{k=1}^K \frac{1}{k} -C\sum_{k=1}^K \frac{1}{k^3}\simeq \log K.
\end{equation*}
The constants $A'$ and $C$ emerge from \eqref{eq:A2,5} and \eqref{eq:A2}. This justifies \eqref{L1_jpol_claim} for $d=1$ and $\al\geq -1/2$.

Let us now assume that $-1/2> \al\geq \be>-1$. We set $\te=\pi/2$ and use \eqref{eq:A5} obtaining
\begin{equation*}
\jtpolk\Big(\frac{\pi}{2}\Big)=2^{\eta} \sqrt{\frac{2}{\pi}}\cos\Big( \frac{k\pi}{2}+\frac{(\be-\al)\pi}{4}\Big)+O(k^{-1}),\qquad  k\in\NN_+.
\end{equation*}
Thus, with $C>0$ resulting from the above asymptotic,
\begin{equation*}
\sum_{k\in\NN_+} \frac{\Big|\jtpolk\big(\frac{\pi}{2}\big)\Big|}{k}\geq 2^{\eta} \sqrt{\frac{2}{\pi}} \sum_{k\in\NN_+} \frac{\Big|\cos\big( \frac{k\pi}{2}+\frac{(\be-\al)\pi}{4}\big)\Big|}{k}-C\sum_{k\in\NN_+} \frac{1}{k^{2}}=\infty.
\end{equation*}

In the multi-dimensional situation we also restrict ourselves to the case $\al_i\geq\be_i$ for $1\leq i\leq d$. 
Let $j\in\{0,1,\ldots,d\}$ be such that $|\{i\colon \al_i<-1/2\}|=j$. Without any loss of generality we assume that $\al_i<-1/2$ for $i=1,\ldots,j$. Observe that claim \eqref{L1_jpol_claim} takes the form
\begin{equation*}
\sup_{\te\in(0,\pi)^d}\sum_{n\in\NN^d_+}\frac{|\jtpol(\te)|}{\ven^{d+(d-j)/2 +\sum_{i=j+1}^d \al_i}}=\infty.
\end{equation*}

For $K$ and $c$ as before we set $\te_K=\big(\underbrace{\frac{\pi}{2},\ldots,\frac{\pi}{2}}_j,\frac{c}{K},\ldots,\frac{c}{K}\big)$. We estimate
\begin{align*}
\sup_{K\in\NN_+}&\sum_{n\in\NN^d_+}\frac{|\jtpol(\te_K)|}{\ven^{d+(d-j)/2 +\sum_{i=j+1}^d \al_i}}\\
&\geq\sup_{K\in\NN} \bigg(A'' \sum_{n\in\{1,\ldots,K \}^d}\frac{\prod_{i=1}^j  \big|\cos\big(\frac{n_i \pi}{2}+\frac{\be_i-\al_i)\pi}{4}\big)\big| \prod_{s=j+1}^d (n_s+1)^{\al_s+1/2}}{\ven^{d+(d-j)/2+\sum_{i=j+1}^d \al_i}}\\
&\hspace{1cm}-C\sum_{i=1}^j \sum_{n\in\{1,\ldots,K \}^d}\frac{1}{n_i \ven^{d+(d-j)/2+\sum_{i=j+1}^d \al_i}}\\
&\hspace{1cm}-C\sum_{i=j+1}^d \sum_{n\in\{1,\ldots,K \}^d}\frac{n_i^{\al_i-3/2}}{\ven^{d+(d-j)/2+\sum_{i=j+1}^d \al_i}}\bigg)\\
&\geq C_{\al,\be}\sum_{n\in\NN_+^d}\frac{\prod_{i=1}^j  \big|\cos\big(\frac{n_i \pi}{2}+\eta_i\big)\big| \prod_{s=j+1}^d (n_s+1)^{\al_s+1/2}}{\ven^{d+(d-j)/2 +\sum_{i=j+1}^d \al_i}}-C\sum_{i=1}^d\sum_{n\in\NN^d_+}\frac{1}{n_i^{1/2}\ven^d},
\end{align*}
where $A'',C>0$ are constants resulting from the applied asymptotics. Divergence of the first series on the right hand side of the latter inequality follows from Lemma \ref{series_divergence_lemma}, whereas the second one is clearly convergent. 

This finishes the justification of \eqref{L1_jpol_claim} and completes the proof of the theorem.
\end{proof}

\section{Jacobi function setting}\label{S4}
Throughout this section we assume $\al,\be\geq -1/2$ (or $\al,\be\in[-1/2,\infty)^d$ in the multi-dimensional case), as the functions $\jfun$ are not in $L^\infty$ if $\min(\al,\be)<-1/2$ (and analogously in higher dimensions). It can be easily checked that in such situation Hardy's inequality and its $L^1$-analogue do not hold. 

We define
\begin{equation*}
\jRop(\te,\va)=\sum_{k=0}^{\infty}r^k \jfunk(\te)\jfunk(\va),\qquad \te,\va\in(0,\pi),\ r\in(0,1).
\end{equation*}
Note that we have the relation
\begin{equation}\label{R&H_kernel_relation}
\jRop(\te,\va)=r^{-\eta} \mathbb{H}^{\al,\be}_{-\log r}(\te,\va). 
\end{equation}

Firstly we will present an auxiliary lemma. Similar result was obtained in \cite[Lemma~1]{KanjinSato}. Nevertheless, our proof seems to be rather shorter, and thus we give it below.

{\lm\label{Jacobi_lemma_dif} If $\al,\be\geq -1/2$, then
\begin{equation*}
|\jfunk(\te)-\jfunk(\te')|\lesssim (k+1)^{\max(\al,\be)+1/2}\left((k+1)|\te-\te'|+|\te-\te'|^{\al+1/2}+|\te-\te'|^{\be+1/2}\right),
\end{equation*}
uniformly in $\te,\te'\in(0,\pi)$ and $k\in\NN$. If $\al=-1/2$ ($\be=-1/2$), then the second (third) component on the right hand side of the estimate can be omitted.}
\begin{proof}
Fix $\al,\be\geq -1/2$. Note that \eqref{jfunk_definition} and \eqref{jtPol_estim} yield
\begin{align*}
&\hspace{-1cm}\big|\jfunk(\te)-\jfunk(\te')\big|\\
&\leq\big(\sin\frac{\te}{2}\big)^{\al+1/2}\big(\cos\frac{\te}{2}\big)^{\be+1/2} \big|\jtpolk(\te)-\jtpolk(\te')\big|\\
&\hspace{0.4cm}+\big(\sin\frac{\te}{2}\big)^{\al+1/2}\Big|\big(\cos\frac{\te}{2}\big)^{\be+1/2}-\big(\cos\frac{\te'}{2}\big)^{\be+1/2}\Big|\big|\jtpolk(\te')\big|\\
&\hspace{0.4cm}+\Big|\big(\sin\frac{\te}{2}\big)^{\al+1/2}-\big(\sin\frac{\te'}{2}\big)^{\al+1/2}\Big|\big(\cos\frac{\te'}{2}\big)^{\be+1/2}\big|\jtpolk(\te')\big|\\
&\lesssim \big|\jtpolk(\te)-\jtpolk(\te')\big|+(k+1)^{\max(\al,\be)+1/2}\Big|\big(\cos\frac{\te}{2}\big)^{\be+1/2}-\big(\cos\frac{\te'}{2}\big)^{\be+1/2}\Big|\\
&\hspace{0.4cm} +(k+1)^{\max(\al,\be)+1/2}\Big|\big(\sin\frac{\te}{2}\big)^{\al+1/2}-\big(\sin\frac{\te'}{2}\big)^{\al+1/2}\Big|.
\end{align*}

Applying the mean value theorem and \eqref{jtpol_differ_estim} we obtain
$$\big|\jtpolk(\te)-\jtpolk(\te')\big|\lesssim |\te-\te'|(k+1)^{\max(\al,\be)+3/2},\qquad k\in\NN,\ \te,\te'\in(0,\pi). $$
Moreover, we remark that for $\eta\geq 0$ we have
\begin{equation}\label{sine_Holder_continuity}
\big|\sin^{\eta} x-\sin^{\eta}y\big|\lesssim |x-y|+|x-y|^\eta\lesssim |x-y|^\eta, \qquad x,y\in[0,\pi/2],
\end{equation}
and the same for the cosines in place of the sines. Indeed, for $\eta\geq 1$ one can use the mean value theorem, whereas for $\eta\in(0,1)$ it suffices to use the $\eta$-H\"older continuity of the function $x\to\sin^\eta x$ on $[0,\pi/2]$.

Combining the above gives the claim.
\end{proof}

{\prop\label{prop_jfun_R_nonsymmetrized} If $\al,\be\geq -1/2$, then 
	\begin{equation}\label{jfun:prop:claim:1}
		\sup_{\te\in(0,\pi)}\Vert \jRop(\te,\cdot)\Vert_{L^2(0,\pi)}\lesssim (1-r)^{-1/2},
	\end{equation}
	uniformly in $r\in(0,1)$, and
	\begin{align}
	&\Vert \jRop(\te,\cdot)-\jRop(\te',\cdot)\Vert_{L^2(0,\pi)}\nonumber\\
	&\qquad\lesssim |\te-\te'|(1-r)^{-3/2}+|\te-\te'|^{\al+1/2}(1-r)^{-(1+\al)}+|\te-\te'|^{\be+1/2}(1-r)^{-(1+\be)}\label{jfun:prop:claim:2},
	\end{align}
	uniformly in $\te,\te'\in(0,\pi)$ and $r\in(0,1)$ (if $\al=-1/2$ or $\be=-1/2$, then we omit respectively the second or the third component from the right hand side of the last estimate).}
\begin{proof}
	Fix $\al,\be\geq -1/2$. We apply Parseval's identity and \eqref{jfun_estim} receiving
	\begin{equation*}
	\sup_{\te\in(0,\pi)}\Vert \jRop(\te,\cdot)\Vert_{L^2((0,\pi))}\lesssim \Big(\sum_{k\in\NN} r^{2k} \Big)^{1/2}\lesssim (1-r)^{-1/2},
	\end{equation*}
	uniformly in $r\in(0,1)$. This finishes the verification of \eqref{jfun:prop:claim:1}.
	
For the proof of \eqref{jfun:prop:claim:2} note that Parseval's identity and Lemma \ref{Jacobi_lemma_dif} yield
\begin{align*}
\Vert \jRop(\te,\cdot)-\jRop(\te',\cdot)\Vert_{L^2(0,\pi)}&\lesssim \Big( \sum_{k=1}^\infty 2^{-2k}|\jfunk(\te)-\jfunk(\te')|^2 \Big)^{1/2}\\
&\lesssim |\te-\te'|+|\te-\te'|^{\al+1/2}+|\te-\te'|^{\be+1/2},
\end{align*}
uniformly in $\te,\te'\in(0,\pi)$ and $r\in(0,1/2]$. 

From now on we assume that $r\in (1/2,1)$. We apply \eqref{R&H_kernel_relation} and calculate assuming $0<\te\leq\te'<\pi$
\begin{align*}
|\jRop(\te,\va)-\jRop(\te',\va)|&\lesssim\Big|\int_{\te}^{\te'} \frac{d}{d\omega} \mathbb{H}^{\al,\be}_{-\log r}(\omega,\va)\,d\omega\Big|\\
&\leq\int_{\te}^{\te'} \Big(\big|D_1^{\al,\be}(\omega,\va)\big|+\big| D_2^{\al,\be}(\omega,\va)\big|+\big|D_3^{\al,\be}(\omega,\va)\big|\Big)\,d\omega,
\end{align*}
where (see \eqref{Jacobi_kernel_relation})
\begin{align*}
D_1^{\al,\be}(\omega,\va)&=\frac{2\al+1}{4}\cos\frac{\omega}{2}\sin\frac{\va}{2}\Big(\sin\frac{\omega}{2}\sin\frac{\va}{2}\Big)^{\al-1/2}\Big(\cos\frac{\omega}{2}\cos\frac{\va}{2}\Big)^{\be+1/2}\mathcal{H}^{\al,\be}_{-\log r}(\omega,\va),\\
D_2^{\al,\be}(\omega,\va)&=\frac{2\be+1}{4}\sin\frac{\omega}{2}\cos\frac{\va}{2}\Big(\sin\frac{\omega}{2}\sin\frac{\va}{2}\Big)^{\al+1/2}\Big(\cos\frac{\omega}{2}\cos\frac{\va}{2}\Big)^{\be-1/2}\mathcal{H}^{\al,\be}_{-\log r}(\omega,\va),\\
D_3^{\al,\be}(\omega,\va)&=\Big(\sin\frac{\omega}{2}\sin\frac{\va}{2}\Big)^{\al+1/2}\Big(\cos\frac{\omega}{2}\cos\frac{\va}{2}\Big)^{\be+1/2}\frac{d}{d\omega}\mathcal{H}^{\al,\be}_{-\log r}(\omega,\va).
\end{align*}

Firstly, we estimate $D_3^{\al,\be}(\omega,\va)$. Applying Theorem \ref{NowakSjogren_J-P_kernel_estim} we get
\begin{align*}
	|D_3^{\al,\be}(\omega,\va)|\lesssim& \Big(\frac{\omega\va}{(1-r)^2+\omega^2+\va^2}\Big)^{\al+1/2}\Big(\frac{(\pi-\omega)(\pi-\va)}{(1-r)^2+(\pi-\omega)^2+(\pi-\va)^2}\Big)^{\be+1/2}\\
	&\times\frac{(1-r)}{((1-r)^2+(\omega-\va)^2)^{3/2}}\\
	\lesssim& (1-r)^{-1}\frac{1}{((1-r)^2+(\omega-\va)^2)^{1/2}}
\end{align*}
Hence, using Minkowski's inequality we obtain
\begin{align*}
	\Big\Vert \int_{\te}^{\te'}D_3^{\al,\be}(\omega,\cdot)\,d\omega\Big\Vert_{L^2(0,\pi)}
	\lesssim& |\te-\te'|(1-r)^{-1} \sup_{\omega\in(0,\pi)}\bigg(\int_0^\pi
	\frac{d\va}{(1-r)^2+(\omega-\va)^2}\bigg)^{1/2}.
\end{align*}
Note that
\begin{equation}\label{integral_arctg}
	\sup_{\omega\in(0,\pi)}\int_0^\pi\frac{d\va}{(1-r)^2+(\omega-\va)^2}\lesssim (1-r)^{-1}.
\end{equation}
This gives the desired estimate for $D_3^{\al,\be}(\omega,\va)$.

%Firstly, we estimate $D_3^{\al,\be}(\omega,\va)$. Note that for $\al,\be\geq -1/2$ we have the estimate \cite[proof of Theorem 2.2]{NowakSjogren}
%\begin{equation*}
%	|\partial_\te \tpois(\te,\va)|\lesssim t^{-1}\tpois(\te,\va),\qquad \te,\va\in(0,\pi),\ t\in(0,1).
%\end{equation*}
%Hence, \eqref{R&H_kernel_relation_pol} and \eqref{Jacobi_kernel_relation} we get
%\begin{align*}
%|D_3^{\al,\be}(\omega,\va)|&\lesssim\Big(\sin\frac{\omega}{2}\sin\frac{\va}{2}\Big)^{\al+1/2}\Big(\cos\frac{\omega}{2}\cos\frac{\va}{2}\Big)^{\be+1/2} (1-r)^{-1} \mathcal{H}^{\al,\be}_{-\log r}(\omega,\va)\\
%&= (1-r)^{-1}\mathbb{H}^{\al,\be}_{-\log r}(\omega,\va).
%\end{align*}
%Hence, by \eqref{R&H_kernel_relation} and Proposition \ref{prop_jfun_R_nonsymmetrized} we see that 
%\begin{align*}
%\Big\Vert \int_{\te}^{\te'} D_3^{\al,\be}(\omega,\cdot)\,d\omega\Big\Vert_{L^2(0,\pi)}&\lesssim (1-r)^{-1}|\te-\te'|\sup_{\omega\in(0,\pi)} \big\Vert R^{\al,\be}_{r}(\omega,\cdot)\big\Vert_{L^2(0,\pi)}\\
%&\lesssim |\te-\te'| (1-r)^{-3/2}.
%\end{align*}

Let us consider $D_1^{\al,\be}$. If $\al=-1/2$, then $D_1^{\al,\be}(\omega,\va)\equiv 0$. Secondly, if $\al\geq 1/2$, then Theorem \ref{NowakSjogren_J-P_kernel_estim} yields
\begin{align*}
&D_1^{\al,\be}(\omega,\va)\\
&\lesssim \frac{(\omega\va)^{\al-1/2}(\pi-\omega)\va}{((1-r)^2+\omega^2+\va^2)^{\al+1/2}}\Big(\frac{(\pi-\omega)(\pi-\va)}{(1-r)^2+(\pi-\omega)^2+(\pi-\va)^2}\Big)^{\be+1/2}\frac{1-r}{(1-r)^2+(\omega-\va)^2}\\
&\lesssim (1-r)^{-1} \frac{(1-r)\va}{(1-r)^2+\va^2}\frac{1}{((1-r)^2+(\omega-\va)^2)^{1/2}}\\
&\lesssim (1-r)^{-1}\frac{1}{((1-r)^2+(\omega-\va)^2)^{1/2}}.
\end{align*}
Thus, invoking Minkowski's inequality and \eqref{integral_arctg} we receive
\begin{align*}
\Big\Vert \int_{\te}^{\te'} D_1^{\al,\be}(\omega,\cdot)\,d\omega\Big\Vert_{L^2(0,\pi)}\lesssim |\te-\te'|(1-r)^{-3/2}.
\end{align*}
Lastly, for $\al\in(-1/2,1/2)$ we also use Theorem \ref{NowakSjogren_J-P_kernel_estim} and obtain
\begin{align*}
D_1^{\al,\be}(\omega,\va)&\lesssim\frac{\va^{\al+1/2}}{((1-r)^2+\omega^2+\va^2)^{\al+1/2}}\frac{1-r}{(1-r)^2+(\omega-\va)^2}\frac{d}{d\omega}\Big(\sin\frac{\omega}{2}\Big)^{\al+1/2}\\
&\lesssim (1-r)^{-\al-1/2}\frac{1-r}{(1-r)^2+(\omega-\va)^2}\frac{d}{d\omega}\Big(\sin\frac{\omega}{2}\Big)^{\al+1/2}.
\end{align*}
Hence, \eqref{sine_Holder_continuity} and \eqref{integral_arctg} imply
\begin{align*}
\Big\Vert& \int_{\te}^{\te'} D_1^{\al,\be}(\omega,\cdot)\,d\omega\Big\Vert_{L^2(0,\pi)}\\
&\lesssim (1-r)^{-(\al+1/2)} \bigg|\Big(\sin\frac{\te}{2}\Big)^{\al+1/2}-\Big(\sin\frac{\te'}{2}\Big)^{\al+1/2}\bigg| \Big(\int_0^\pi \frac{d\va}{(1-r)^2+\va^2}\Big)^{1/2}\\
&\lesssim |\te-\te'|^{\al+1/2} (1-r)^{-(\al+1)} .
\end{align*}

The computations for $D_2^{\al,\be}(\omega,\va)$ are similar to those for $D_1^{\al,\be}(\omega,\va)$. Combining the obtained results justifies \eqref{jfun:prop:claim:2}. This completes the proof of the proposition.
\end{proof}

{\thm\label{Jacobi_fun_nonsym_mainthm} Let $\al,\be\in[-1/2,\infty)^d$. Hardy's inequality for the Jacobi trigonometric functions holds with the admissible exponent $E=d$, namely 
\begin{equation*}
\sum_{n\in\NN^d} \frac{|\langle f,\jfun\rangle|}{(\ven+1)^d}\lesssim \Vert f\Vert_{H^1((0,\pi)^d)},
\end{equation*}
uniformly in $f\in H^1((0,\pi)^d)$. The admissible exponent is sharp, i.e. for any $\ve>0$ there exists $f\in H^1((0,\pi)^d)$ such that
\begin{equation*}
\sum_{n\in\NN^d} \frac{|\langle f,\jfun\rangle|}{(\ven+1)^{d-\ve}}=\infty.
\end{equation*} }

\begin{proof}
Fix $\al,\be\in[ -1/2,\infty)^d$. In order to prove the inequality we shall again invoke Theorem \ref{general_thm}. The proper assumptions are satisfied. Indeed, Proposition \ref{prop_jfun_R_nonsymmetrized} gives for
\begin{equation*}
\jRop(\te,\va)=\sum_{n\in\NN^d} r^{\ven}\jfun(\te)\jfun(\va),\qquad \te,\va\in(0,\pi)^d,\ r\in(0,1),
\end{equation*}
the following estimate
\begin{align*}
&\Vert \jRop(\te,\cdot)-\jRop(\te',\cdot)\Vert_{L^2((0,\pi)^d)}\\
&\lesssim |\te-\te'|(1-r)^{-(d+2)/2}\\
&\hspace{1cm}+\sum_{i=1}^d \Big( |\te-\te'|^{\al_i+1/2}(1-r)^{-\al_i-(d+1)/2}+|\te-\te'|^{\be_i+1/2}(1-r)^{-\be_i-(d+1)/2}\Big),
\end{align*} 
where we omit the components corresponding to $\al_i=-1/2$ and $\be_i=-1/2$. Thus, the corresponding parameter $\gamma$ (see \eqref{parameter_gamma}) is equal to $(d+2)/2$ and 
$$\Delta=\{1,\al_1+1/2,\ldots,\al_d+1/2,\be_1+1/2,\ldots,\be_d+1/2\}, $$
where, again, we exclude the zeros. Furthermore, $N=d$  (see \eqref{parameter_N}) for Lebesgue measure, and thus by \eqref{E_formula} $E=d$, which finishes the verification of Hardy's inequality. 

Now we shall justify sharpness  (for $d=1$; if $d\geq 1$ it suffices to take the atom defined by the tensor product of the one-dimensional ones). We may assume that $\al\geq \be$. In the light of Theorem \ref{Jacobi_pol_nonsym_mainthm} and the fact that $\phi_k^{-1/2,-1/2}\equiv \mathcal{P}_k^{-1/2,-1/2}$, we can restrict ourselves to the case $\al>-1/2$.

Similarly as in the proof of Theorem \ref{Jacobi_pol_nonsym_mainthm}, it suffices, for a given $K\in\NN_+$, to construct an $H^1((0,\pi))$-atom $a_K$ such that
\begin{equation*}
	\sum_{k=1}^K \frac{|\langle a_K,\jfunk\rangle|}{k^{1-\ve}}\gtrsim K^\ve,
\end{equation*}
with the implicit constant independent of $K$.

Let us fix large $K\in\NN_+$. We define for $\delta\in(0,1/2)$
\begin{equation*}
a_K(\te)=\left\{ \begin{array}{rl}
-K/c,&\te\in(0,c\delta K^{-1}],\\
\frac{\delta}{1-\delta}K/c, &\te\in(c\delta K^{-1},cK^{-1}],\\
0,&\text{otherwise},
\end{array}
\right.
\end{equation*}
where $c$ is the constant from \eqref{eq:A2}. It can be easily checked that $a_K$ is an $H^1((0,\pi))$-atom. 

%Note that \eqref{eq:A1} implies that
%\[\jfunk(\te)=((k+\eta)\te)^{1/2}J_{\al}((k+\eta)\te)+R_k^{\al,\be}(\te),\qquad 0<\te<cK^{-1}, \]
%where $|R_k^{\al,\be}(\te)|\leq C \te^{\al+1/2}k^{\al-3/2}$. 
We apply \eqref{eq:A1} along with \eqref{eq:A2} and compute for $K/2\leq k\leq K$
\begin{align*}
\int_0^{cK^{-1}}a_K(\te)\jfunk(\te)\,d\te&\geq -\frac{K}{Ac} \int_0^{c\delta K^{-1}}(\te k)^{\al+1/2}\,d\te\\
&\qquad+\frac{AK}{c}\frac{\delta}{1-\delta}\int_{c\delta K^{-1}}^{cK^{-1}}(\te k)^{\al+1/2}\,d\te\\
&\qquad-Ck^{\al-3/2} K^{-\al-1/2} \Vert a\Vert_{L^1((0,\pi))} \\
&\geq \Big(\frac{k}{K}\Big)^{\al+1/2}\bigg(\frac{2c^{\al+1/2}\delta }{(2\al+3)A}\Big( \frac{A^2(1-\delta^{\al+3/2})}{(1-\delta)}-\delta^{\al+1/2}\Big)-\frac{4C}{K^2}\bigg),
\end{align*}
where $A>0$ is the constant from \eqref{eq:A2} and $C>0$ results from the asymptotic in \eqref{eq:A1}.

Thus, for sufficiently small $\delta$ and $K$ large enough, we obtain
\begin{equation*}
|\langle a_K,\jfunk\rangle|\geq \Big(\frac{k}{K}\Big)^{\al+1/2}\Big(\frac{A(1-\delta^{\al+3/2}) c^{\al+1/2}\delta}{(2\al+3)(1-\delta)}-\frac{4C}{K^2}\Big)\gtrsim 1,\qquad K/2\leq k\leq K.
\end{equation*}
In the first inequality we used the assumption that $\al> -1/2$. Hence,
\begin{equation*}
\sum_{k=1}^K\frac{|\langle a_K,\jfunk\rangle|}{k^{1-\ve}}\gtrsim \sum_{k=\lceil K/2\rceil}^K \frac{1}{k^{1-\ve}}\gtrsim K^\ve.
\end{equation*}

This finishes the proof of the theorem.
%This proves \eqref{H1_jfun_claim} for $\al>-1/2$. If $\al=-1/2$ and $\be>-1/2$ we can conduct a similar reasoning using the symmetry $\jfunk(\te)=\phi_k^{\be,\al}(\pi-\te)$. Therefore, it remains to deals with the case $\al=\be=-1/2$.

%Recall that $\phi^{-1/2,-1/2}_k$ are the Chebyshev polynomials given by (cf. \cite[(4.1.7)]{Szego})
%\begin{equation*}
%\phi_k^{-1/2,-1/2}(\te)=\left\{ \begin{array}{ll}
%\sqrt{2/\pi},&k=0,\\
%\sqrt{1/\pi}\cos k\te,&k\geq 1.
%\end{array}\right.
%\end{equation*}
%Hence, for $k=1,\ldots,K$ we estimate invoking the mean value theorem
%\begin{align*}
%\Big| \int_0^{\frac{c}{K}}a(\te)\frac{1}{\sqrt{\pi}}\cos k\te\, d\te\Big|&=\Big| \int_0^{\frac{c}{K}}a(\te)\frac{1}{\sqrt{\pi}}\big(\cos k\te-\cos\frac{kc\delta}{N}\big)\, d\te\Big|\\
%&=\frac{Kk}{c\sqrt{\pi}}\Big| \int_0^{\frac{c\delta}{K}}\big(\te-\frac{c\delta}{K}\big)\sin k\xi_\te\, d\te-\frac{\delta}{1-\delta}\int_{\frac{c\delta}{K}}^{\frac{c}{K}}\big(\te-\frac{c\delta}{K}\big)\sin k\xi_\te\, d\te\Big|\\
%&\gtrsim k^2\delta\int_{\frac{c\delta}{K}}^{\frac{c}{K}}\big(\te-\frac{c\delta}{K}\big)\, d\te\\
%&=\frac{k^2\delta^2 c^2}{2K^2(1-\delta)}(1-\delta^2-2\delta+2\delta^2)\\
%&\gtrsim k^2K^{-2},
%\end{align*}
%where $\xi_\te$ lies between $\te$ and $\frac{c\delta}{K}$, and we used the fact that $\sin k\xi_\te\gtrsim \frac{kc\delta}{K}$ for $\frac{c\delta}{K}\leq \te\leq \frac{c}{K}$. Thus,
%\begin{equation*}
%\sum_{k=1}^K \frac{|\langle a,\phi_k^{-1/2,-1/2}\rangle|}{k^{1-\ve}}\gtrsim K^{-2}\sum_{k=1}^K k^{1+\ve}\gtrsim K^{\ve}.
%\end{equation*}

\end{proof}

{\thm\label{nonsym_jfun_L1thm} Let $\al,\be\in[-1/2,\infty)^d$ and fix $\ve>0$. The following inequality holds
\begin{equation*}
\sum_{n\in\NN^d} \frac{|\langle f,\jfun\rangle|}{(\ven+1)^{d+\ve}}\lesssim \Vert f\Vert_{L^1((0,\pi)^d)},
\end{equation*}
uniformly in $f\in L^1((0,\pi)^d)$. Moreover, the admissible exponent is sharp.}

\begin{proof}
Fix $\al,\be\in [-1/2,\infty)^d$. The inequality follows from \eqref{jfun_estim}.

In order to justify sharpness it is sufficient to verify that (again see \cite[Lemma~1]{Kanjin2})
\begin{equation}\label{L1_jfun_claim}
\sup_{\te\in(0,\pi)^d}\sum_{n\in\NN^d}\frac{|\jfun(\te)|}{(\ven+1)^d}=\infty.
\end{equation}
We use the multi-dimensional version of \eqref{eq:A6} for $\te=(\pi/2,\ldots,\pi/2)$, namely
\begin{equation*}
\jfun\Big(\frac{\pi}{2},\ldots,\frac{\pi}{2}\Big)=\Big(\frac{2}{\pi}\Big)^{d/2} \prod_{i=1}^d \cos\Big(\frac{n_i\pi}{2}+\frac{(\be_i-\al_i)\pi}{4}\Big)+\sum_{i=1}^d O(n_i^{-1}),\qquad n\in\NN_+^d.
\end{equation*}
Thus, with $C>0$ emerging from the above asymptotic, we simply obtain
\begin{equation*}
\sum_{n\in\NN^d}\frac{\big|\jfun\big(\frac{\pi}{2}\big)\big|}{(\ven+1)^d}\geq \Big(\frac{2}{\pi}\Big)^{d/2}\sum_{n\in\NN_+^d}\frac{\prod_{i=1}^d \Big|\cos\big(\frac{n_i\pi}{2}+\frac{(\be_i-\al_i)\pi}{4}\big)\Big|}{\ven^d} -C\sum_{i=1}^d \sum_{n\in\NN_+^d}\frac{1}{\ven^d n_i}.
\end{equation*}
Note that the reminder converges absolutely. Moreover, Lemma \ref{series_divergence_lemma} implies that the first series on the right hand side of the above inequality diverges.

This finishes justification of \eqref{L1_jfun_claim} and hence the proof of the theorem.

\end{proof}

\section{Symmetrized Jacobi Settings}\label{S5}
In this section we shall establish Hardy's inequality and its $L^1$-analogue for the symmetrized Jacobi trigonometric polynomials and functions. In fact, we will justify that the appropriate claims follow from the corresponding ones from the previous sections.

{\thm\label{sym_Jacobi_fun_mainthm} Let $\al,\be\in[-1/2,\infty)^d$. Sharp Hardy's inequality for the symmetrized Jacobi trigonometric functions holds with the admissible exponent $E=d$. Its sharp $L^1$-analogue with $E=d+\ve$, $\ve>0$, is also valid.}
\begin{proof}
For a function $f$ defined on $(-\pi,\pi)^d$ we introduce the following decomposition:
\begin{equation*}
f=\sum_{\sigma\in\{0,1\}^d} f_\sigma,\qquad f_\sigma =2^{-d} \sum_{\epsilon\in\{-1,1\}^d} \epsilon^\sigma f(\epsilon\cdot),
\end{equation*}
where $\epsilon^\sigma=\prod_{i=1}^d \epsilon_i^{\sigma_i}$. Note that if $\sigma_i=0$, then $f_\sigma$ is an even function with respect to the $i$-th variable, and if $\sigma_i=1$, then $f_\sigma$ is odd with respect to the $i$-th variable.

Let $\mf\colon\NN^d\to\{0,1\}^d$ be the function returning the parity of an integer multi-index, namely $\mf(n)=n\mod 2$. For any $\sigma\in\{0,1\}^d$ and $\al\in(-1,\infty)^d$ we write $\al_\sigma=((\al_1)_{\sigma_1},\ldots,(\al_d)_{\sigma_d})$, where
\begin{equation*}
(\al_i)_{\sigma_i}=\left\{\begin{array}{ll}
\al_i,& \sigma_i=0,\\
\al_i+1,& \sigma_i=1, 
\end{array}\right.
\end{equation*}
and similarly for $\be$. 

We remark that 
\begin{equation*}
|\langle f,\sjfun\rangle|=|\langle f_{\mf(n)},\sjfun\rangle| \simeq |\langle f_{\mf(n)}^+,(\sjfun)^+\rangle_+|\simeq |\langle f_{\mf(n)}^+,\phi_{\lfloor n/2\rfloor}^{\al_{\mf(n)},\be_{\mf(n)}}\rangle_+|,
\end{equation*}
where $f^+$ denotes the restriction of $f$ to $(0,\pi)^d$, the inner product $\langle\cdot,\cdot\rangle$ is taken in $L^2((-\pi,\pi)^d)$, and $\langle\cdot,\cdot\rangle_+$ is taken in $L^2((0,\pi)^d)$. Here $\lfloor n/2\rfloor = (\lfloor n_1/2\rfloor, \ldots,\lfloor n_d/2\rfloor)$. Hence, applying Theorem \ref{Jacobi_fun_nonsym_mainthm} we receive
\begin{align*}
\sum_{n\in\NN^d} \frac{|\langle f,\sjfun\rangle|}{(\ven+1)^d}\simeq \sum_{\sigma\in\{0,1\}^d} \sum_{\mf(n)=\sigma} \frac{|\langle f_\sigma^+,\phi_{\lfloor n/2\rfloor}^{\al_\sigma,\be_\sigma}\rangle_+|}{(\ven+1)^d}\lesssim \sum_{\sigma\in\{0,1\}^d} \Vert f_\sigma^+\Vert_{H^1((0,\pi)^d)}\lesssim \Vert f\Vert_{H^1((-\pi,\pi)^d)},
\end{align*}
where in the last estimate we used the inequality (cf. \cite[(5)]{Plewa3})
\begin{equation*}
\Vert f_\sigma^+\Vert_{H^1((0,\pi)^d)}\leq \Vert f\Vert_{H^1((-\pi,\pi)^d)}.
\end{equation*}

Sharpness in this setting follows easily from sharpness for the Jacobi trigonometric functions (Theorem \ref{Jacobi_fun_nonsym_mainthm}). Moreover, the $L^1$-analogue can be justified very similarly, therefore we skip this part of the proof. 
\end{proof}

A similar theorem holds for the symmetrized Jacobi trigonometric polynomials. However, it is significantly harder to deduce this from the corresponding results for the Jacobi polynomials, as the admissible exponents for $\jtpol$ and $\mathcal{P}_{n}^{\al+1,\be+1}$ are different. Moreover, the underlying spaces are not the same and this produces some difficulties.

{\thm\label{sym_Jacobi_pol_mainthm} If $\al,\be\in(-1,\infty)^d$, then sharp Hardy's inequality for the symmetrized Jacobi trigonometric polynomials holds with the admissible exponent 
$$E=3d/2 +\sum_{i=1}^d \max(\al_i,\be_i,-1/2).$$
If we add $\ve>0$ to this exponent, then also the $L^1$-inequality holds.}
\begin{proof}
One can immediately deduce sharpness in both claims from the results in Section \ref{S3}. Hence, it suffices to justify Hardy's inequality (the $L^1$-analogue follows). 

Observe that if we will proceed similarly as in the proof of Theorem \ref{sym_Jacobi_fun_mainthm}, we will encounter the trigonometric polynomials
\begin{equation*}
Q_k^{\al,\be}(\te)=\frac{\sin\te}{2}\mathcal{P}_k^{\al+1,\be+1}(\te),\qquad \te\in (0,\pi).
\end{equation*}
Note that $\{Q_k^{\al,\be}\}_{k\in\NN}$ forms an orthonormal basis in $L^2((0,\pi),\mu_{\al,\be})$. Moreover, if we set 
\begin{equation*}
\tilde{R}_r^{\al,\be}(\te,\va)=\sum_{k=0}^\infty r^k Q_k^{\al,\be}(\te)Q_k^{\al,\be}(\va),\qquad \te,\va\in(0,\pi),\ r\in(0,1),
\end{equation*}
then it is easy to justify that the results analogous to Proposition \ref{prop_jpol_R_nonsymmetrized} hold. Indeed, we have
\begin{equation*}
	|Q_k^{\al,\be}(\te)|\lesssim (k+1)^{1/2+\max(\al,\be,-1/2)},\qquad |\partial_\te Q_k^{\al,\be}(\te)|\lesssim (k+1)^{3/2+\max(\al,\be,-1/2)},
\end{equation*}
uniformly in $\te\in(0,\pi)$ and $k\in\NN$. Thus, Parseval's identity and \eqref{series_estim_eq} yield
\begin{equation*}
\sup_{\te\in(0,\pi)} \Vert \tilde{R}_r^{\al,\be}(\te,\cdot)\Vert_{L^2((0,\pi),\,\mu_{\al,\be})}\lesssim (1-r)^{-1-\max(\al,\be,-1/2)},
\end{equation*}
\begin{equation*}
\sup_{\te\in(0,\pi)} \Vert \partial_\te\tilde{R}_r^{\al,\be}(\te,\cdot)\Vert_{L^2((0,\pi),\,\mu_{\al,\be})}\lesssim (1-r)^{-2-\max(\al,\be,-1/2)},
\end{equation*}
uniformly in $r\in(0,1)$.

 Hence, Hardy's inequalities associated with the Jacobi trigonometric polynomials and the polynomials $\{Q_k^{\al,\be}\}_{k\in\NN}$ hold with the same exponent and, consequently, so does Hardy's inequality associated with the symmetrized Jacobi trigonometric polynomials. 
 
 This concludes the proof.
\end{proof}

\appendix
\section{Asymptotic estimates of the Hilb and Darboux types}
\numberwithin{equation}{section}
Fix $\al,\be\in(-1,\infty)$. We will make use of the following version of Hilb's type formula (see \cite[(5.7)]{StempakTohoku'02} and \cite[8.21.17]{Szego} for the original version)
\begin{equation}\label{eq:A1}
	\jfunk(\te)=((k+\eta)\te)^{1/2}J_{\al}((k+\eta)\te)+O(\te^{\al+1/2}k^{\al-3/2}),\quad \te\in(0,ck^{-1}),\, k\in\NN_+,
\end{equation}
where $c$ is a small positive fixed constant, which is independent of $k$, $\te$, and $J_\al$ denotes the Bessel function of first kind (see \cite[(1.71.1)]{Szego}). Recall that $\eta=(\al+\be+1)/2$. We remark that we do not use the final form of this formula \cite[(4.6)]{StempakTohoku'02}, since the remainder is too large. 

Using the known fact that $J_\al(z)\simeq z^{\al}$ for sufficiently small $z$, we obtain
\begin{equation}\label{eq:A2}
	A (k\te)^{\al+1/2}\leq ((k+\eta)\te)^{1/2}J_{\al}((k+\eta)\te)\leq \frac{1}{A} (k\te)^{\al+1/2},\quad \te\in(0,ck^{-1}),\ k\in\NN_+,
\end{equation}
with $c>0$ possibly smaller than in \eqref{eq:A1} and some constant $A>0$.

We will also need a version of \eqref{eq:A1} for the Jacobi trigonometric polynomials. By \eqref{jfunk_definition} we can simply write
\begin{equation}\label{eq:A2,5}
\jtpolk(\te)=\Big(\sin\frac{\te}{2}\Big)^{-\al-1/2}\Big(\cos\frac{\te}{2}\Big)^{-\be-1/2}((k+\eta)\te)^{1/2}J_{\al}((k+\eta)\te)+O(k^{\al-3/2}),
\end{equation}
uniformly in $\te\in(0,c/k)$, $k\in\NN_+$.

Applying \eqref{jtPol_differentiation_formula} we get
\begin{align}
-\frac{d}{d\te}\jtpolk(\te)=&\sqrt{k(k+2\eta)}\Big(\sin\frac{\te}{2}\Big)^{-\al-1/2}\Big(\cos\frac{\te}{2}\Big)^{-\be-1/2}((k+\eta)\te)^{1/2}\nonumber\\
&\times J_{\al+1}((k+\eta)\te)+O(\te k^{\al+1/2}),\qquad \te\in(0,ck^{-1}),\ k\in\NN_+. \label{eq:A3}
\end{align}

Analogously to \eqref{eq:A2} we have
\begin{align}
	B \te k^{\al+5/2}&\leq \sqrt{k(k+2\eta)}\Big(\sin\frac{\te}{2}\Big)^{-\al-1/2}\nonumber\\
	&\hspace{1cm}\times \Big(\cos\frac{\te}{2}\Big)^{-\be-1/2}\sqrt{(k+\eta)\te}J_{\al+1}((k+\eta)\te)\leq \frac{1}{B} \te k^{\al+5/2},\label{eq:A4}
\end{align}
uniformly in $\te\in (0,ck^{-1})$ and $k\in\NN_+$, for some $B>0$ depending on $\al$ and $\be$.

Now we pass to the Darboux formula (cf. \cite[(8.21.18)]{Szego}) rewritten to the form
\begin{align*}
\jtpolk(\te)&=c_k^{\al,\be}k^{-1/2}\pi^{-1/2}\Big(\sin\frac{\te}{2}\Big)^{-\al-1/2}\Big(\cos\frac{\te}{2}\Big)^{-\be-1/2}\cos\Big( (k+\eta)\te-\frac{(2\al+1)\pi}{4}\Big)\\
&\hspace{1cm}+O(k^{-1}),
\end{align*}
uniformly for $k\in\NN_+$ and $\te$ separated from $0$ and $\pi$, say $\te\in[\pi/6,5\pi/6]$. Here $c_k^{\al,\be}$ is the normalizing constant given in \eqref{normalizing_constant}. Observe that using Stirling's approximation we get
\[\big|c_k^{\al,\be}k^{-1/2}-\sqrt{2}\big|\lesssim k^{-1},\qquad k\in\NN_+. \]
Hence, we may write
\begin{equation}\label{eq:A5}
\jtpolk(\te)=\sqrt{\frac{2}{\pi}}\Big(\sin\frac{\te}{2}\Big)^{-\al-1/2}\Big(\cos\frac{\te}{2}\Big)^{-\be-1/2}\cos\Big( (k+\eta)\te-\frac{(2\al+1)\pi}{4}\Big)+O(k^{-1}),
\end{equation}
uniformly in $\te\in[\pi/6,5\pi/6]$ and $k\in\NN_+$.

Similar asymptotic estimate holds for the Jacobi trigonometric functions. Indeed, \eqref{jfunk_definition} yields
\begin{equation}\label{eq:A6}
\jfunk(\te)=\sqrt{\frac{2}{\pi}}\cos\Big( (k+\eta)\te-\frac{(2\al+1)\pi}{4}\Big)+O(k^{-1}),\qquad \te\in\Big[\frac{\pi}{6},\frac{5\pi}{6}\Big].
\end{equation}


\begin{thebibliography}{999}
	
	\bibitem{BalasRadha}
	R. Balasubramanian, R. Radha,
	{\it Hardy-type inequalities for Hermite expansions},
	J. Inequal. Pure Appl. Math. 6 (2005), 1--4.
	
	\bibitem{CastroNowakSzarek}
	A. J. Castro, A. Nowak, T. Z. Szarek,
	{\it Riesz--Jacobi Transforms as Principal Value Integrals},
	J. Fourier Anal. Appl. 22 (2016), 493--541.
	

	\bibitem{CoifmanWeiss}
	R. R. Coifman, G. Weiss,
	{\it Extension of Hardy spaces and their use in analysis},
	Bull. Amer. Math Soc. 83 (1977), 569--645.

	\bibitem{HardyLittlewood}
	G. H. Hardy, J. E. Littlewood,
	{\it Some new properties of Fourier constants},
	Math. Annal. 97 (1927), 159--209.

	\bibitem{Kanjin1}
	Y. Kanjin,
	{\it Hardy's inequalities for Hermite and Laguerre expansions},
	Bull. London Math. Soc. 29 (1997), 331--337.
	
	\bibitem{Kanjin2}
	Y. Kanjin,
	{\it Hardy's inequalities for Hermite and Laguerre expansions revisited},
	J. Math. Soc. Japan 63 (2011), 753--767.
	
%	\bibitem{KanjinSato2}
%	Y. Kanjin, K. Sato,
%	{\it Paley's inequality the for Jacobi expansions},
%	Bull. London Math. Soc. 33 (2001), 483--491.
	
	\bibitem{KanjinSato}
	Y. Kanjin, K. Sato,
	{\it Hardy's inequality for Jacobi expansions},
	Math. Inequal. Appl. 7 (2004), 551--555.

	\bibitem{Langowski}
	B. Langowski,
	{\it Harmonic analysis operators related to symmetrized Jacobi expansions for all admissible parameters},
	Acta Math. Hungar. 150 (2016), 49--82.
	
	\bibitem{Lebedev}
	N. N. Lebedev,
	{\it Special functions and their applications}, Dover, New York, 1972.
	
	\bibitem{LiYuShi}
	Z. Li, Y. Yu, Y. Shi,
	{\it The Hardy inequality for Hermite expansions},
	J. Fourier Anal. Appl. 21 (2015), 267--280.	

	\bibitem{Muckenhoupt2}
	B. Muckenhoupt,
	{\it Transplantation theorems and multiplier theorems for Jacobi series},
	Mem. Am. Math. Soc. 64, 1986.
	
	\bibitem{NowakSjogren}
	A. Nowak, P. Sj\"{o}gren,
	{\it Calder\'o{}n-Zygmund Operators Related to Jacobi Expansions},
	J. Fourier Anal. Appl. 18 (2012), 717--749.
	
	\bibitem{NowakSjogren2}
	A. Nowak, P. Sj\"{o}gren,
	{\it Sharp estimates of the Jacobi heat kernel},
	Stud. Math. 218 (2013), 219--244.
	
	\bibitem{NowakSjogrenSzarek}
	A. Nowak, P. Sj\"{o}gren, T. Z. Szarek,
	{\it Analysis Related to All Admissible Type Parameters in the Jacobi Setting},
	Constr. Approx. 41 (2015), 185--218.

	\bibitem{NowakStempak(symm)}
	A. Nowak, K. Stempak,
	{\it A symmetrized conjugacy scheme for orthogonal expansions},
	Proc. Roy. Soc. Edinburgh Sect. A 143 (2013), 427--443.

	\bibitem{Plewa}
	P. Plewa,
	{\it Hardy's type inequality for Laguerre expansions of Hermite type},
	J. Fourier Anal. Appl. (2018), DOI: 10.1007/s00041-018-9642-2.

	\bibitem{Plewa2}
	P. Plewa,
	{\it Sharp Hardy's type inequality for Laguerre expansions},
	preprint (2018), arXiv:1810.08138.

	\bibitem{Plewa3}
	P. Plewa,
	{\it On Hardy's inequality for Hermite expansions},
	Taiwanese J. Math. (2019), DOI: 10.11650/tjm/190601.
	
	\bibitem{Radha}
	R. Radha,
	{\it Hardy type inequalities},
	Taiwanese J. Math. 4 (2000), 447--456.
	
	\bibitem{RadhaThangavelu}
	R. Radha, S. Thangavelu
	{\it Hardy's inequalities for Hermite and Laguerre expansions},
	Proc. Amer. Math. Soc. 132 (2004), 3525--3536.
	
	\bibitem{Satake}
	M. Satake,
	{\it Hardy's inequalities for Laguerre expansions},
	J. Math. Soc. Japan 52 (2000), 17--24.
	
	\bibitem{StempakTohoku'02}
	K. Stempak,
	{\it On connections between Hankel, Laguerre and Jacobi Transplantations},
	Tohoku Math. J. 54 (2002), 471--493.
	
	\bibitem{Szego}
	G. Szeg\"{o},
	{\it Orthogonal polynomials},
	Amer. Math. Soc. Colloq. Publ. Vol 23, fourth edition, Providence, 1975.
	
	\bibitem{Thangavelu}
	S. Thangavelu,
	{\it Hermite and Laguerre Expansions},
	Princeton University Press, Princeton (1993).
	


\end{thebibliography}
\end{document}